\documentclass[12pt,draftcls,onecolumn]{IEEEtran}

\usepackage{graphicx}
\usepackage{url}
 \usepackage[pdftex]{color}
\usepackage{amsmath}
\usepackage{amsfonts} 
\usepackage{array}
\usepackage{cancel}
\newtheorem{teor}{Theorem}[section]
\newtheorem{corr}{Corollary}[section]
\newtheorem{propo}{Proposition}[section]
\newtheorem{lemm}{Lemma}[section]

\newtheorem{rem}{Remark}
\usepackage{soul}
\usepackage{tabularx}
\usepackage{multirow}

\newcommand{\tp}{^{\top}}

\newcommand{\beq}{\begin{equation}}
\newcommand{\eeq}{\end{equation}}
\newcommand{\bea}{\begin{eqnarray}}
\newcommand{\eea}{\end{eqnarray}}

\newcommand{\bsea}{\begin{subeqnarray}}
\newcommand{\esea}{\end{subeqnarray}}
\newcommand{\nn}{\nonumber}

\newcommand{\proof}{\noindent {\it Proof. }}

\newcommand{\qed}{\hfill $\Box$ \vskip 2ex}

\def\bmat{\left[ \begin{array}}
\def\emat{\end{array} \right]}
\DeclareMathOperator{\tr}{tr} 
\DeclareMathOperator{\dd}{diag^{2}} 
\DeclareMathOperator{\ofd}{ofd}
\DeclareMathOperator{\diag}{diag} 
\DeclareMathOperator*{\argmin}{arg\,min}
\definecolor{Royalblue}{cmyk}{1,0.30,0.2,0.2}
\newcommand{\agu}{\color{black}}
\newcommand{\vale}{\color{black}}
\newcommand{\alg}[1]{\begin{align}#1\end{align}}
\newcommand{\mz}{\color{black}} 
 
\newcommand{\gus}{\color{black}}

\begin{document}

\title{ Factor Models with Real Data:\\ a Robust Estimation of the Number of Factors}

\author{Valentina~Ciccone, Augusto~Ferrante, Mattia~Zorzi \thanks{V. Ciccone, A. Ferrante, and M. Zorzi are with   the Department of Information
Engineering, University of Padova, Padova, Italy; e-mail: {\tt\small valentina.ciccone@dei.unipd.it} (V. Ciccone); {\tt\small augusto@dei.unipd.it} (A. Ferrante); {\tt\small zorzimat@dei.unipd.it} (M. Zorzi).}}

\markboth{DRAFT}{Shell \MakeLowercase{\textit{et al.}}: Bare Demo of IEEEtran.cls for Journals}

\maketitle

\begin{abstract}
Factor models are a very efficient way to describe high dimensional vectors of data in terms of a small number of common relevant factors. This problem, which is of fundamental importance in many disciplines, is usually reformulated in mathematical terms as  follows. We are given the covariance matrix $\Sigma$ of the available data. $\Sigma$ must be additively decomposed as the sum of two positive semidefinite matrices $D$ and $L$: $D$ | that accounts for the  idiosyncratic noise affecting the knowledge of each component of the available vector of data | must be diagonal and $L$ must have the smallest possible rank  in order to describe the available data in terms of the smallest possible number of independent factors.

In practice, however, the matrix $\Sigma$ is never {\agu known and therefore it} must be estimated from the data so that only an {\agu approximation} of 
$\Sigma$ is actually  available. This paper discusses the {\mz issues} that arise from this uncertainty and provides a strategy to deal {\mz with the problem of {\agu robustly} estimating the number of factors.}
 \end{abstract}

\begin{IEEEkeywords}
factor analysis; nuclear norm; convex optimization; duality theory.
\end{IEEEkeywords}

\section{Introduction}

Describing a large amount of data by means of a small number of factors carrying most of the information is {\vale an important problem in} modern data analysis with applications ranging in all fields of science. One of the classical methods for this purpose is to resort to {\em factor models} that were first developed at the beginning of the last century  by Spearman \cite{spearman_1904} in the framework of the so-called {\em mental tests} as an attempt at ``the procedure of eliciting verifiable facts''  in determining  psychical tendencies from the tests results.
From this first seed a rich stream of literature was developed at the interface between psychology and mathematics with the main focus on the case of a single common factor underlying the available data:  necessary and sufficient conditions for the data to be compatible with a single common factor were derived in \cite{BURT_1909,Spearman-Holzinger-24}, see also \cite{Bekker-deLeeuw_1987} and references therein for a detailed historical reconstruction of the derivation of these conditions. The interest for this kind of model has grown rapidly also outside the psychology community and analysis of factor models, or {\em factor analysis} has become an important tool in statistics, econometrics, systems theory and many engineering fields \cite{KALMAN_SELECTION_ECONOMETRICS_1983,SCHUPPEN_1986,Bekker-deLeeuw_1987,PICCI_1987,GEWEKE_DYNAMIC,PICCI_PINZONI_1986,Pena_BOX__1987,
ON_THE_PENA_BOX_MODEL_2001,deistler2007,DEISTLER_1997,ANDERSON_DEISTLER_1993,
SARGENT_SIM_1977,forni_lippi_2001,ENGLE_ONE_FACTOR_1981,WATSON_ALGORITHMS_1983}, {\vale \cite{mclachlan1997wiley}}; see also the more recent papers \cite{Bottegal-Picci,MFA,DEISTLER_2015,fan2013large}, {\vale \cite{bertsimas2017certifiably,delgado2014rank}} where many other references are listed.  A detailed geometric description of this problem is presented  in \cite{scherrer1998structure}.
In the seminal paper \cite{anderson1956statistical} a maximum likelihood approach in a statistical testing framework is proposed.

In the original formulation the construction of a factor model is equivalent to the mathematical problem of additively decomposing a given positive definite  matrix $\Sigma$ | modeling the covariance of the data | as 
\beq\label{add-dec}
\Sigma=L+D
\eeq
 where both $L$ and $D$ are positive semidefinite, and $D$ | modeling the covariance of the idiosyncratic noise | is diagonal.
The rank of $L$ is the number of (latent or hidden) common factors that explain the available data.
One of the key aspects of factor analysis is to determine the minimum number of latent factors  or, equivalently, a decomposition 
(\ref{add-dec}) where the rank of $L$ is minimal. 
This  is therefore a particular case of a matrix additive decomposition problem that arises naturally in numerous frameworks and have therefore received a great deal of attention,
see \cite{Chandrasekaran-Sanghavi-Parrilo-Willsky,agarwal2012noisy,LATENTG,BSL} and references therein.
We hasten to remark that the problem of minimizing the rank of $L$ in the decomposition (\ref{add-dec}) is extremely hard so that, the convex relaxation is usually considered where, in place of the rank, the nuclear norm (i.e. the trace) of $L$ is minimized. This is a very good approximation that most often returns, with reasonable computational burden, a solution $L$ with minimum rank.

\begin{figure}[htb]
	\centering
	\includegraphics[trim={3.5cm 2.5cm 2cm 1.8cm},clip,width=0.5\textwidth]{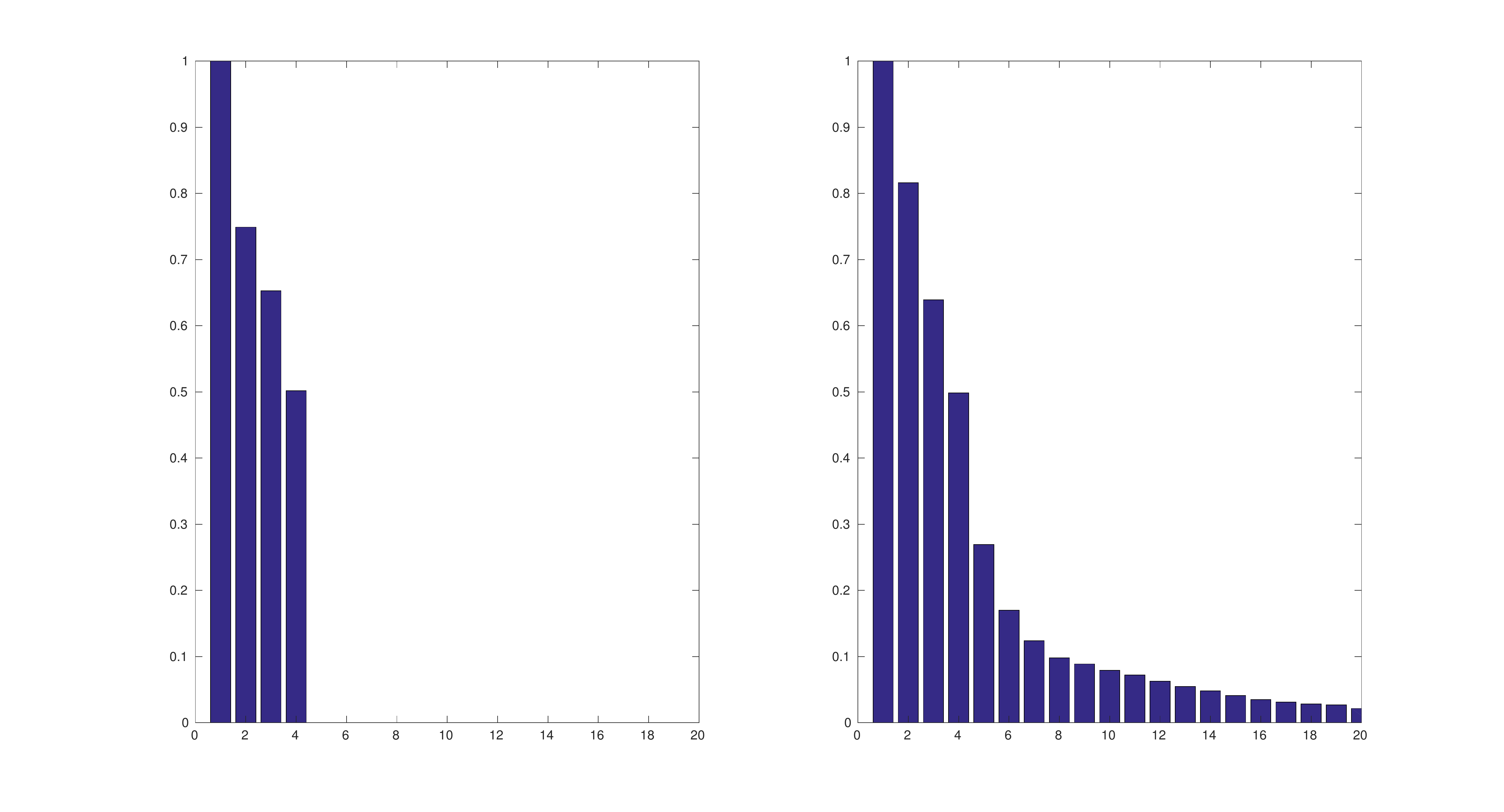}
	\caption{ {\vale First $20$ singular values of the  matrices  $L$  (on the left) and $\hat{L}$ (on the right)    obtained by applying the minimum trace factor analysis decomposition algorithm to  a ``true'' covariance matrix $\Sigma\in\mathbb{R}^{40\times 40}$
of a model with $r=4$ latent factors and to an estimate $\hat{\Sigma}$ of $\Sigma$
obtained by generating $N=1000$ independent samples from a normal distribution 
$\mathcal{N}(0,\Sigma)$ and computing the corresponding sample covariance, respectively. Notice that the recovered matrix $L$ using as input the true $\Sigma$ is, up to negligible numerical errors, equal to the true low rank matrix.}}
\label{fig:barplot1}
\end{figure}

In \cite{THURSTONE_1935},\cite{Kelley_1928},\cite{Ledermann_1937} an upper bound $r(n)$ | known as {\em Ledermann bound} | was proposed for the minimal rank  $r_m(\Sigma)$ of $L$ 
in terms if the dimension $n$ of the matrix  $\Sigma$:
$$
r_m(\Sigma)\leq r(n):= \Big{\lfloor} \frac{2n+1-\sqrt{8n+1}}{2} \Big{\rfloor}.
$$
This bound, however, is based on heuristics that  have never been proven rigorously; a {\em p\'etale de rose} is the  prize for a positive demonstration of this fact \cite{Hakim1976}.\footnote{Indeed, not only is a rigorous proof  missing but a precise statement is also needed. In fact, some further assumptions must be added for the validity of this bound as counterexamples can, otherwise, be easily produced \cite{Guttman1958}.}
Interestingly, almost half a century later in \cite{RANK_REDUCIBILITY_1982} a related result was established:  the set of matrices $\Sigma$ for which $r_m(\Sigma)< r(n)$ has zero Lebesgue measure. 
{\agu As a consequence of this result 
we have the following observation that may be regarded as the basic premise of our effort. When $n$ is large the Ledermann bound $r(n)$ is not much smaller than $n$. \gus Therefore, even if our data do come from a factor model with a small number $r$ of latent factors,  only a set of zero measure of $\hat{\Sigma}$ in a neighbourhood of $\Sigma$  can be decomposed in such a way that the corresponding $L$ matrix in its decomposition (\ref{add-dec}) has rank $r$. Thus, unless we know $\Sigma$ with absolute precision,
we cannot rely only on the decomposition (\ref{add-dec}) to recover such $r$.}
{\vale An example of this phenomenon is illustrated in Figure \ref{fig:barplot1}.} 

{\vale The problem of estimating $r$ from an estimate $\hat{\Sigma}$ of $\Sigma$ is therefore of crucial importance and has been addressed in \cite{bai2002determining} and \cite{lam2012factor} by means of statistical methods.
{\vale A similar issue has been {\agu addressed also} in \cite{ning2015linear} {\agu in the framework of the robustness of Frisch scheme.}}
Here, we propose an {\agu alternative optimization-based} approach which {\agu is based only on the estimate $\hat{\Sigma}$ {\agu and}  takes {\mz into account} the uncertainty of this estimate}.} {\agu Hence, even if we can start from $N$ $n$-dimensional vectors (observations)
the data of our problem are just the sample covariance  $\hat{\Sigma}$ of these vectors and their number $N$. These
two quantities summarize all the relevant information for our method in which}  we  {\em compute} the matrix $\Sigma$
in such a way that the {\agu trace} of $L$ in its additive decomposition (\ref{add-dec})  is minimized under a constraint limiting the Kullback-Leibler divergence between $\Sigma$ and $\hat{\Sigma}$ to a prescribed tolerance that depends on the precision of our estimate $\hat{\Sigma}$  and hence may be reliably chosen on the basis of the data numerosity $N$.

The {\vale proposed} problem is analyzed by resorting {\agu to}  duality theory. The {\agu dual analysis} is delicate to carry over, but yields a problem whose solution can be efficiently computed by {\agu employing}  an alternating direction method of multipliers (ADMM) algorithm. Moreover, the dual problem provides a necessary and sufficient condition for the uniqueness of the solution of the original problem.

The paper is organized as follows. In the Section \ref{sec_pb_form} we recall the classical approach to factor analysis and, from it, we derive the formulation of our factor analysis  problem. In Section \ref{sec_choice_delta} we describe how to establish, for a desired tolerance, an upper bound on the aforementioned Kullback-Leibler divergence. In Section \ref{sec_dual} we derive a dual formulation of our problem. In Section \ref{sec_ex} we prove existence and uniqueness of the solution for the dual problem. Then, in Section \ref{sec_recovery} we show how to recover the solution of the primal problem. 
In Section \ref{algorithm} we present the numerical algorithm for solving the dual problem, while in Section \ref{simulations} the results of numerical simulations and an application to a real world example are presented.
Finally, some conclusions are provided.
The less instructive proofs that are essentially based on calculations are deferred to the Appendix.

Some of the results of this paper have been presented in preliminary form and mostly without proof in our conference paper \cite{NF_CDC}. 

{\em Notation:} Given a vector space $\cal V$ and a subspace ${\cal W}\subset {\cal V}$, we denote by
${\cal W}^\bot$ the orthogonal complement of ${\cal W}$ in ${\cal V}$.
Given a matrix $M$, we denote its transpose by $M\tp$; if $M$ is a square matrix 
$\tr(M)$ denotes its trace, i.e. the sum of the elements in the main diagonal of $M$;
moreover, $|M|$ denotes the determinant of $M$ and $\sigma(M)$ denotes the spectrum of $M$, i.e. the set of its eigenvalues. {  We denote the spectral norm of $M$ as $\Vert M \Vert_2$.}
We endow the space of square real matrices with the following inner product: for $A,B\in{\mathbb R}^{n\times n}$, $\langle A,B\rangle:= \tr(A^TB)$. 
The kernel of a matrix (or of a linear operator) is denoted by $\ker(\cdot)$.
The symbol $\mathbf{Q}_n$ denotes the vector space of real symmetric matrices of size $n$. If $X\in\mathbf{Q}_n$ is positive definite or positive semi-definite we write $X\succ 0$ or $X\succeq 0$, respectively.
Moreover, we denote by $\mathbf{D}_n$ the vector space of  diagonal matrices of size $n$; $\mathbf{D}_n$ is clearly a subspace of $\mathbf{Q}_n$ and we denote by
$\mathbf{M}_n:=\mathbf{D}_n^{\bot}$ the orthogonal complement of $\mathbf{D}_n$ in $\mathbf{Q}_n$ (with respect to the inner product just defined).
It is easy to see that 
$\mathbf{M}_n$ is the vector space of symmetric matrices of size $n$ having all the elements on the main diagonal equal to zero. We denote by $\diag(\cdot)$ both the operator mapping $n$ real elements $d_i, i=1,...,n$ into the diagonal matrix having the $d_i$'s as elements in its main diagonal and the operator mapping a matrix $M\in\mathbb{R}^{n\times n}$ into an $n$-dimensional vector containing the diagonal elements of $M$. Then $\diag\diag(\cdot)$, that we denote by
 $\dd(\cdot)$, is the (orthogonal projection) operator mapping  a square matrix $M$ into a diagonal matrix of the same size having the same main diagonal of $M$.
We denote by $\ofd(\cdot)$ the self-adjoint operator orthogonally projecting $\mathbf{Q}_n$ onto $\mathbf{M}_n$, i.e. if $M\in  \mathbf{Q}_n$, $\ofd(M)$ is the matrix of $\mathbf{M}_n$ in which each off-diagonal element is equal to the corresponding element of $M$ (and each diagonal element is clearly zero). Finally, we denote by $\otimes$ the Kronecker product between two matrices   and by vec(X) the vectorization of a matrix X formed by stacking the columns of X into a single column vector.

\section{Problem Formulation}\label{sec_pb_form}

We consider a standard factor model in its static linear formulation
 \alg{\label{fact_model} y&=Ax+z}
 where $A\in \mathbb{R}^{n\times r}$, with $r<<n$, is the factor loading matrix, $x$ represents the (independent) latent factors and $z$ is the idiosyncratic component. $x$ and $z$  are independent Gaussian random vectors with zero mean and covariance matrix equal to the identity matrix of dimension $r$ and $D\in \mathbf{D}_n$, respectively. Note that, $Ax$ represents the latent variable. 
Consequently, $y$ is a Gaussian random vector with zero mean; we denote by $\Sigma$ its covariance matrix. Since $x$ and $z$ are independent we get that {\agu $\Sigma$ may be additively decomposed as in \eqref{add-dec},} where $L:=AA^\top$ and $D$ are the covariance matrices of $Ax$ and $z$, respectively. Thus, $L$  has rank equal to $r$, and $D$ is diagonal.

The objective of factor analysis consists in finding {\vale the most parsimonious ``low-rank plus diagonal'' decomposition  of $\Sigma$, that is {\agu a decomposition \eqref{add-dec} } for which the rank of $L$ is minimal}. This amounts to solve the minimum rank problem
\begin{equation}
\label{rank_min}
\begin{aligned}
\min_{L,D\in \mathbf{Q}_n} \quad & \text{ rank}(L)\\
\text{subject to }\quad & L,D\succeq 0\\
					& D \in \mathbf{D}_n\\
					& {\Sigma}=L+D
\end{aligned}
\end{equation}
which is, however, a hard problem. 
A well-known {\vale and widely used heuristic is the} convex relaxation of \eqref{rank_min}{\vale, i.e.} the trace minimization problem
 \begin{equation}
\label{original_problem}
\begin{aligned}
\min_{L,D\in \mathbf{Q}_n} \quad & \text{ tr}(L)\\
\text{subject to }\quad & L,D\succeq 0\\
					    & D\in \mathbf{D}_n\\
					    & {\Sigma}=L+D.
\end{aligned}
\end{equation}
{\mz The substitution of the rank with the trace} is justified by the fact that $\tr(L)$, i.e. the nuclear norm of $L$, is the convex hull of $\text{ rank}(L)$ over the set ${\cal S}:=\{ L\in \mathbf{Q}_n \; \mathrm{s.t.} \; \| L\|_2\leq 1\}$, \cite{FAZEL_MIN_RANK_APPLICATIONS_2002}. {\vale The relation between Problem \eqref{rank_min} and Problem \eqref{original_problem} has been first studied in {\gus\cite{MIN_RANK_SHAPIRO_1982} and
while these two problems are, in general, not equivalent, very often they have the same solution.} }

In practice, however, matrix $\Sigma$ is not known and needs to be estimated from a $N$-length realization (i.e. a data record) $\mathrm{y}_1 \ldots \mathrm{y}_N$ of $y$. The typical choice is to take the sample covariance estimate \alg{\hat \Sigma:=\frac{1}{N}\sum_{k=1}^N  \mathrm{y}_k  \mathrm{y}_k\tp \label{samplecov}}
which is statistically consistent, i.e. the corresponding estimator almost surely converges to $\Sigma$ as $N$ tends to infinity. {\vale As discussed in the Introduction, by replacing $\Sigma$ with $\hat \Sigma$ 
 the solution,  in terms of minimum rank, will rapidly degrade. 
Indeed a delicate problem in factor analysis is the one of estimating the number of factors. Such a problem has been addressed by several important contributions, see the seminal works of Bai and Ng \cite{bai2002determining} and of Lam and Yao \cite{lam2012factor} and the references therein.
Our objective is to address the same problem from a different perspective. In fact, we propose an optimization problem whose solution provides an estimate of the minimum number of factors by introducing an appropriate model for the error in the estimation of $\Sigma$. 
} 
{\agu This model is based on an auxiliary Gaussian random vector $\hat y$  with zero mean and covariance matrix $\hat \Sigma$
that is regarded as a ``model approximation'' for $y$. To account for {\vale the estimation} uncertainty, we assume that
the distribution of $y$ (that is completely specified by its covariance matrix and hence is referred to by $\Sigma$)  belongs to a ``ball'' centred in  $\hat y$
\alg{{\cal B}:=\{ \Sigma\in \mathbf{Q}_n \; \mathrm{s.t.}\; \Sigma\succ 0,\;{\cal D}_{KL}(\Sigma\|\hat \Sigma)\leq \delta /2\}} which is formed by placing a bound (i.e. tolerance) on the {\em Kullback-Leibler} divergence between $y$ and $\hat y$:
\alg{ {\cal D}_{KL}(\Sigma\|\hat \Sigma):= \frac{1}{2}\left(-\log |\Sigma|+\log |\hat \Sigma|+\tr(\Sigma \hat \Sigma^{-1})-n\right).\nn}
}
This way to deal with model uncertainty has been successfully {\vale applied} in econometrics for model mispecification \cite{ROBUSTNESS_HANSENSARGENT_2008} and in robust filtering \cite{ROBUST_FILTERING,ROBUST_KALMAN_2017,LEVY_NIKOUKHAH_2004,ROBUST_CONV_2015,OPTIMALITY_ZORZI,convtau}. Accordingly, {\vale in order to estimate the minimum number of factors, we propose the following ``robustification'' of the minimum trace problem:}
\begin{equation}
\label{primal_before_changing_variables}
\begin{aligned}
\min_{\Sigma,L,D\in  \mathbf{Q}_n} & \text{ tr}(L)\\
\text{subject to  } &  L,D\succeq 0\\
					& D\in \mathbf{D}_n\\
					& \Sigma=L+D\\
					& \Sigma\in \cal B.
\end{aligned}
\end{equation}

Note that, in (\ref{primal_before_changing_variables}) we can eliminate variable $D$, obtaining the equivalent problem 
\begin{equation}
\label{primal}
\begin{aligned}
\min_{L,\Sigma\in  \mathbf{Q}_n} \quad & \text{ tr}(L)\\
\text{subject to  } \quad & L,\Sigma-L\succeq 0 \\
					& \ofd(\Sigma-L)=0 \\
					& \Sigma\succ 0\\
					&2 \mathcal{D}_{KL}(\Sigma||\hat{\Sigma})\leq \delta. 
\end{aligned}
\end{equation}
{\mz It is worth noting that an alternative to Problem (\ref{primal}) is to consider $\mathcal{D}_{KL}(\Sigma||\hat{\Sigma})$ as {\agu a} penalty term in the objective function rather than as a constraint. Such approach, however, would require a cross validation procedure to set the regularization parameter $\lambda$, i.e. we {\agu would} have to solve an optimization problem for many values of $\lambda$. In contrast, the proposed problem is solved only once provided that $\delta$ is chosen in a suitable way, see the next section.}

\section{The choice of $\delta$}\label{sec_choice_delta}
 The tolerance $\delta$ may be chosen by taking into account the accuracy of the estimate $\hat \Sigma$ of $\Sigma$ which, in turn, depends on the numerosity of the available data. This can be done by choosing a probability $\alpha\in (0,1)$ and a neighborhood of ``radius'' $\delta_\alpha$ (in the Kullback-Leibler topology) centered in $\hat \Sigma$ containing the ``true'' $\Sigma$ with probability $\alpha$.
 The Kullback-Leibler divergence in \eqref{primal} is a function of the estimated sample covariance and as such its accuracy depends crucially on the numerosity of the available data.  To asses this accuracy
 we propose an  approach that hinges on  the following scale-invariance property of the  Kullback-Leibler divergence.
\begin{lemm}
Let $y_i\sim\mathcal{N}(0,\Sigma),\, i=1,...,N$ be i.i.d. {\em \mz random variables} {\agu taking values in} $\mathbb{R}^{n}$ and define the sample covariance estimator
\[
 \hat{{\mathbf \Sigma}}=\dfrac{1}{N}\sum_{i=1}^{N}y_i y_i^{T}.
\]
The Kullback-Leibler divergence between  $\Sigma$ and $\hat{{\mathbf \Sigma}}$
is a random variable whose distribution depends only of the number $N$ of random variables  and on the dimension $n$ of each random variable.
\end{lemm}
\proof
We have
\[\hat{{\mathbf \Sigma}}=\dfrac{1}{N}\sum_{i=1}^{N}y_i y_i^{T}=\dfrac{1}{N}\Sigma^{1/2}\sum_{i=1}^{N}\tilde y_i \tilde y_i^{T}\Sigma^{1/2}=\Sigma^{1/2}Q_N\Sigma^{1/2}\]
with $\tilde y_i=\Sigma^{-1/2} y_i\sim\mathcal{N}(0,I_n)$ and $Q_N:=\dfrac{1}{N}\sum_{i=1}^{N}\tilde y_i \tilde y_i^{T},$  is a random matrix taking values in $\mathbf{Q}_n$. {\mz Notice that at this point $\tilde y$ are normalized Gaussian random vectors and hence do not depend on the data nor on $\Sigma$. Thus, $Q_N$ is a random matrix whose distribution only depends on $N$ and $n$ (see Section \ref{sec:G_ens} for more details)}. Hence, the Kullback-Leibler divergence between $\Sigma$ and the sample covariance estimator is 
\begin{equation}
\begin{aligned}
\label{D_1}
d:&=\mathcal{D}_{KL}(\Sigma\Vert\hat{{\mathbf \Sigma}})= \dfrac{1}{2}\big(\log(|\hat{{\mathbf \Sigma}}\Sigma^{-1}|)+\text{tr}(\Sigma\hat{{\mathbf \Sigma}}^{-1}) -n\big) \\
& = \dfrac{1}{2}\big(\log\left(|Q_N|\right) +\text{tr}(Q_N^{-1})-n\big).
\end{aligned}
\end{equation}
\qed
In view of this result we can easily approximate  the distribution of 
the random variable $2d=2\mathcal{D}_{KL}(\Sigma\Vert\hat{{\mathbf \Sigma}})$ 
by a standard Monte Carlo method. 
In particular, we can reliably estimate with arbitrary precision  the value of $\delta$ for which
$Pr\big(2 \mathcal{D}_{KL}(\Sigma||\hat{\Sigma})\leq \delta \big)=\alpha$.
As an alternative to this empiric approach for determining $\delta_{\alpha}$, we can also resort to an analytic one as discussed below.

\subsection{Gaussian Orthogonal Ensemble}\label{sec:G_ens}
Let us focus on the random matrix $Q_N$ that we have defined as $Q_N:=\dfrac{1}{N}\sum_{i=1}^{N} \tilde y_i \tilde y_i^{T}$ where $\tilde y_i\in\mathbb{R}^{n},\, i=1,...,N$, are i.i.d.  random variables distributed as $\mathcal{N}(0,I_n)$.
We now introduce a new matrix $\tilde{Q}_N:= \sqrt{N} (Q_N-I_n)=\sqrt{N}\dfrac{1}{N}\sum_{i=1}^{N} C_i$, where  $C_i:=\tilde y_i \tilde y_i^{T}-I$ are i.i.d. symmetric random matrices with zero mean.
It is immediate to check that for each $i$:  any two distinct elements $[C_i]_{h,j}$ and $[C_i]_{k,l}$ of $C_i$ are uncorrelated as long as they do not occupy symmetric positions, i.e. whenever $(h,j)\neq (l,k)$, and ${\rm Var}\left[[C_i]_{h,j}\right]=\begin{cases}
		1, & \text{if  } h\neq j \\ 2, & \text{if  } h=j 
		\end{cases}$. By the multivariate Central Limit Theorem, we have that  $\tilde{Q}_N$ converges in distribution to the random matrix 
\[ X= \left( \begin{array}{cccc}
\sqrt{2}\xi_{1,1} & \cdots & \cdots & \xi_{1,n} \\
\vdots & \ddots & \, & \vdots \\
\vdots & \, & \ddots & \vdots \\
\xi_{1,n} & \cdots & \cdots & \sqrt{2}\xi_{n,n} \end{array} \right) \in \mathbf{Q}_n, \] 
where   $\lbrace\xi_{i,j}\rbrace$ are i.i.d. Gaussian random variables with mean 0 and variance 1. The set of these matrices is known as the {\em Gaussian Orthogonal Ensemble}, see  \cite{anderson2010introduction}.
It is well known that the joint distribution of the eigenvalues  $\lambda_1(X)\leq...\leq\lambda_n(X)$ of such matrices takes the following form:
\[p(\lambda_1,...,\lambda_n)=\bar{C}_n|\Delta(\lambda)|\prod_{i=1}^{n}\mathrm{e}^{-\lambda_i^{2}/4}\]
where $\lambda:=(\lambda_1,...,\lambda_n)$, $|\Delta(\lambda)|$ is the Vandermonde determinant associated with $\lambda$, which is given by:
\[ |\Delta(\lambda)|=\prod_{i<j}(\lambda_j-\lambda_i)\] 
and $\bar{C}_n$ is defined as:
\[\bar{C}_n=\left(\int_{-\infty}^{+\infty}...\int_{-\infty}^{+\infty}|\Delta(\lambda)|\prod_{i=1}^{n}\mathrm{e}^{-\lambda_i^{2}/4} d\lambda_i \right)^{-1}.\]

It is not difficult to see that  \eqref{D_1} can be rewritten as:
\begin{equation}
d= d(\lambda_1,...,\lambda_n)=\sum_{i=1}^{n} \dfrac{1}{2}\left(\log\left(\dfrac{\lambda_i}{\sqrt{N}}+1\right)-\dfrac{\lambda_i}{\lambda_i+\sqrt{N}} \right)
\end{equation}
where $\lambda_i\in\sigma(\tilde Q_N )$. Then, for a desired $\alpha$, we are interested in finding $\delta_\alpha$ such that 
\[Pr( 2 d\leq \delta_\alpha)=  \alpha.\]
Such a value for $\delta_\alpha$ is given by the cumulative distribution function $F(\cdot)$:
\begin{equation}
\label{CDF}
\begin{aligned}
F(\delta_\alpha)&=Pr( 2d\leq \delta_\alpha)\\
&=\int_{I(\delta_\alpha)} 2d(\lambda_1,...,\lambda_n)p(\lambda_1,...,\lambda_n)d\lambda
\end{aligned}
\end{equation}
where $p(\lambda_1,...,\lambda_n)$ denotes the joint probability density function of the eigenvalues $\lambda_1,...,\lambda_n$ and $I(\delta_\alpha):=\lbrace (\lambda_1,...,\lambda_n): d(\lambda_1,...,\lambda_n)\leq \delta_\alpha /2\rbrace$.
Given $\alpha$ the integral in \eqref{CDF} can be solved numerically for $\delta_\alpha$.

\subsection{An upper bound for $\delta_\alpha$}
If  the chosen level $\alpha$ is too large with respect to the sample size $N$, the computed $\delta_\alpha$ become excessively large so that there are diagonal matrices 
$\Sigma_D$ such that  $2 \mathcal{D}_{KL}(\Sigma_D||\hat{\Sigma})\leq \delta_\alpha$.
In this case Problem (\ref{primal}) admits the trivial solution $L=0$ and $D=\Sigma_D$.
In order to rule out this trivial situation we need to require that the maximum value for $\delta$ in \eqref{primal} is strictly less than a certain $\delta_{max}$ that can be determined as follows: since the trivial solution $L= 0$ would imply a diagonal $\Sigma$, that is $\Sigma=\Sigma_D:=\text{diag}(d_1,...,d_n)>0$, $\delta_{max}$ can be determined by solving the following minimization problem
\begin{align}
\label{min_diagonal}
\delta_{max}:=\min_{\Sigma_D\in\mathbf{D}_n  }  2\mathcal{D}_{KL}(\Sigma_D \Vert \hat{\Sigma}).\end{align}
The following Proposition, whose proof is in Appendix, shows how to solve this problem.
\begin{propo}\label{prop_deltamax} Let $\gamma_i$ denote the $i$-th element in the main diagonal of the inverse of the sample covariance $\hat{\Sigma}^{-1}$. Then, the optimal $\Sigma_D$ which solves the minimization problem in \eqref{min_diagonal} is given by
\[\Sigma_D= \text{diag}(\gamma_1^{-1}, ...,\gamma_n^{-1}).\]
Moreover, $\delta_{max}$ can be determined as 
\begin{equation}
\label{delta_MAX}
\delta_{max}= 2\mathcal{D}_{KL}(\Sigma_D \Vert \hat{\Sigma})= \log|\dd(\hat{\Sigma}^{-1})\hat{\Sigma}|.
\end{equation} \end{propo}

In what follows, we always assume that 
$\delta$ in \eqref{primal} strictly less than $\delta_{max}$, so that the trivial solution $L= 0$ is ruled out.

\section{Dual Problem}\label{sec_dual}
By duality theory, we reformulate the constrained minimization problem in \eqref{primal} as an unconstrained minimization problem.
The associated Lagrangian is
\begin{equation}
\label{lagrangian}
\begin{aligned}
\mathcal{L} (L,& \Sigma, \lambda, \Lambda, \Gamma, \Theta ) \\
= & \text{tr}(L) +\lambda(-\log|\Sigma|+\log|\hat{\Sigma}|-n +\text{tr}(\hat{\Sigma}^{-1}\Sigma)-\delta)\\
& -\text{tr}(\Lambda L)-\text{tr}(\Gamma(\Sigma-L))+\text{tr}{(\Theta \ofd(\Sigma-L)})\\
= & \text{tr}(L) +\lambda(-\log|\Sigma|+\log|\hat{\Sigma}|-n +\text{tr}(\hat{\Sigma}^{-1}\Sigma)-\delta)\\
& -\text{tr}(\Lambda L)-\text{tr}(\Gamma(\Sigma-L))+\text{tr}{(\ofd^{*}(\Theta)(\Sigma-L)})\\
= & \text{tr}(L) +\lambda(-\log|\Sigma|+\log|\hat{\Sigma}|-n +\text{tr}(\hat{\Sigma}^{-1}\Sigma)-\delta)\\
& -\text{tr}(\Lambda L)-\text{tr}(\Gamma(\Sigma-L))+\text{tr}{(\ofd(\Theta)(\Sigma-L)})\\
\end{aligned}
\end{equation}
with $\lambda\in\mathbb{R}, \lambda\geq 0$, and $\Lambda,\Gamma, \Theta \in \mathbf{Q}_n$ with $\Lambda,\Gamma\succeq 0$. In the last equality, we exploited the fact that the operator $\ofd(\cdot)$ is self-adjoint.  
Note that the Lagrangian \eqref{lagrangian} does not include the constraint $\Sigma\succ 0$: as we will see this condition is automatically met by the solution of the dual problem.\\
{\vale Notice also that in \eqref{lagrangian} we can recognize the fit term $2\lambda \mathcal{D}(\Sigma \Vert \hat{\Sigma})=\lambda(-\log|\Sigma|+\log|\hat{\Sigma}|-n +\text{tr}(\hat{\Sigma}^{-1}\Sigma))$ and the term $\text{tr}((I-\Lambda)L+ (\ofd(\Theta)-\Gamma)(\Sigma-L))$ accounting for the complexity in the class of models \eqref{add-dec} which induces low-rank on matrix $L$. Thus \eqref{lagrangian} can be interpreted as an alternative {\mz to the} likelihood function with a complexity term}.

The dual function is defined as the infimum of $\mathcal{L}(L,\Sigma, \lambda, \Lambda, \Gamma, \Theta )$ over $L$ and $\Sigma$.\\
Thanks to the convexity of the Lagrangian, we rely on standard variational methods to characterize the minimum.\\
The first variation of the Lagrangian \eqref{lagrangian} at $\Sigma$ in direction $\delta\Sigma\in\mathbf{Q}_n$ is 
\begin{equation*}
\delta\mathcal{L}(\Sigma;\delta\Sigma) = \text{tr}(-\lambda\Sigma^{-1}\delta\Sigma+\lambda\hat{\Sigma}^{-1}\delta\Sigma-\Gamma\delta\Sigma+\ofd(\Theta)\delta\Sigma).
\end{equation*}
We impose the optimality condition 
\[\delta\mathcal{L}(\Sigma;\delta\Sigma)=0, \qquad\forall\delta\Sigma\in\mathbf{Q}_n,\]
 which is equivalent to require $\text{tr}(-\lambda\Sigma^{-1}\delta\Sigma+\lambda\hat{\Sigma}^{-1}\delta\Sigma-\Gamma\delta\Sigma+\ofd(\Theta)\delta\Sigma) =0$ for all $\delta\Sigma\in\mathbf{Q}_n$, obtaining
\begin{equation}
\label{der_var_wrt_sigma_opt}
\begin{aligned} 
\Sigma=\lambda(\lambda\hat{\Sigma}^{-1}-\Gamma+\ofd(\Theta))&^{-1}
\end{aligned}
\end{equation}
provided that $\lambda\hat{\Sigma}^{-1}-\Gamma+\ofd(\Theta)\succ 0$ and $\lambda>0$, which is clearly equivalent to require that the optimal $\Sigma$ that minimizes the Lagrangian satisfies the constraint $\Sigma\succ 0$.\\

The first variation of the Lagrangian \eqref{lagrangian} at $L$ in direction $\delta L\in\mathbf{Q}_n$ is
\begin{equation*}
\delta\mathcal{L}(L;\delta L)= \text{tr}(\delta L - \Lambda\delta L+\Gamma\delta L -\ofd(\Theta)\delta L).
\end{equation*}
Again, we impose the optimality condition
\[\delta\mathcal{L}(L;\delta L)=0, \qquad \forall\delta L\in\mathbf{Q}_n,\]
which is equivalent to require $\text{tr}(\delta L - \Lambda\delta L+\Gamma\delta L -\ofd(\Theta)\delta L)= 0$ for all $\delta L\in\mathbf{Q}_n$ and we get that 
 \begin{equation}
\label{der_var_wrt_r}
\begin{aligned}
I -\Lambda +\Gamma -\ofd(\Theta)&=0 .
\end{aligned}
\end{equation}

The following result, whose proof is in Appendix, provides a precise formulation of the dual problem.
\begin{propo} \label{prop_dual}The dual problem of (\ref{primal}) is 
\begin{equation}
\label{max_dual}
\max_{(\lambda,\Gamma, \Theta) \in \mathcal{C}_0} J(\lambda,\Gamma,\Theta) 
\end{equation}
where 
\begin{equation*}
\label{dual}
\begin{aligned}
J(\lambda, \Gamma,\Theta):=  &\lambda(\log|(\hat{\Sigma}^{-1}+\lambda^{-1}(\ofd(\Theta)-\Gamma))|  \\
& +\log|\hat{\Sigma}|-\delta) 
\end{aligned}
\end{equation*}
and $\mathcal{C}_0$ is defined as 
\alg{
\label{C0}
\mathcal{C}_0:=\lbrace (\lambda,\Gamma,\Theta):\ &\lambda>0, \ I+\Gamma-\ofd(\Theta)\succeq0, \   
\Gamma\succeq 0,\nn\\ &\hat{\Sigma}^{-1}+\lambda^{-1}(\ofd(\Theta)-\Gamma)\succ 0\rbrace.
}
\end{propo}

\section{Existence and uniqueness of the solution for the dual problem} \label{sec_ex}

We reformulate the maximization problem in \eqref{max_dual} as a minimization problem:
\begin{equation}
\label{Dual}
\begin{aligned}
\min_{(\lambda,\Gamma, \Theta) \in \mathcal{C}_0} \tilde{J}(\lambda,\Gamma, \Theta )
\end{aligned}
\end{equation}
where 
\begin{equation*}
\begin{aligned}
\tilde{J}(\lambda,\Gamma, \Theta)=& \lambda  (-\log|\hat{\Sigma}^{-1}+\lambda^{-1}(\ofd(\Theta)-\Gamma)|\\
&-\log|\hat{\Sigma}|+\delta).
\end{aligned}
\end{equation*}

\subsection{Existence}
As it is often the case, existence of the optimal solution is a very delicate issue.
Our strategy in order to deal with that is to prove that the dual problem in \eqref{Dual} admits  solution.  In doing that we show that we can restrict our set $\mathcal{C}_0$ to a compact set $\mathcal{C}$ over which the minimization problem is equivalent to the one in \eqref{Dual}. Since the objective function is continuous over $\mathcal{C}_0$, and hence over $\mathcal{C}$, by Weierstrass's theorem $\tilde{J}$ admits a minimum.

First, we recall that the operator $\ofd(\cdot)$ is self-adjoint. Moreover, we notice that $\ofd(\cdot)$ is not injective on $\Theta$, thus we can restrict the domain of $\ofd(\cdot)$ to those $\Theta$ such that $\ofd(\cdot)$ is injective. Since $\ofd$ is self-adjoint we have that:
\[\ker(\ofd)=[\text{range }(\ofd)]^{\perp}.\]
Thus, by restricting $\Theta$ to range($\ofd$)$=[\ker(\ofd)]^{\perp}=\mathbf{M}_n$, the map becomes injective. Therefore, without loss of generality, from now on we can safely assume that $\Theta\in\mathbf{M}_n$ so that
$\ofd(\Theta)=\Theta$ and we 
restrict our set $\mathcal{C}_0$ to $\mathcal{C}_{1}$:
\begin{align*}
\mathcal{C}_{1}: =& \lbrace (\lambda,\Gamma,\Theta)\in \mathcal{C}_0: \Theta \in \mathbf{M}_n \rbrace \\
= & \lbrace (\lambda,\Gamma,\Theta):\lambda>0, I+\Gamma-\Theta\succeq 0, \Gamma\succeq 0, \\
& \Theta \in \mathbf{M}_n,  (\hat{\Sigma}^{-1}+\lambda^{-1}(\Theta-\Gamma))\succ 0 \rbrace.
\end{align*}
Moreover, since $\Theta$ and $\Gamma$ enter into the problem always through their difference they cannot be univocally determined individually. However, their difference does. This allows us to restrict $\Gamma$ to the space of the diagonal positive semi-definite matrices.
Indeed, for any sequence $(\lambda_k, \Gamma_k, \Theta_k)_{k\in\mathbb{N}}\in \mathcal{C}_1$ such that $\inf \tilde{J}(\lambda, \Gamma, \Theta)= \lim_{k\rightarrow \infty}\tilde{J}(\lambda_k, \Gamma_k, \Theta_k)$ we can always consider a different sequence $({\lambda}_k, \tilde{\Gamma}_k, \tilde{\Theta}_k)_{k\in\mathbb{N}}$ with $\tilde{\Gamma}_k:= \dd(\Gamma_k)$ and $\tilde{\Theta}_k:=\Theta_k-\ofd(\Gamma_k)$. It is now immediate to check that the new sequence still belongs  to $\mathcal{C}_1$ and that we still have $\inf \tilde{J}(\lambda, \Gamma, \Theta)= \lim_{k\rightarrow \infty}\tilde{J}(\tilde{\lambda}_k, \tilde{\Gamma}_k, \tilde{\Theta}_k)$.
For this reason, we can further restrict our set $\mathcal{C}_1$ to $\mathcal{C}_{2}$:
\begin{align*}
\mathcal{C}_{2}:= & \lbrace (\lambda, \Gamma, \Theta): \lambda>0, I+\Gamma-\Theta\succeq 0, \Gamma\succeq 0, \Gamma\in\mathbf{D}_n,\\
& \Theta \in \mathbf{M}_n, \hat{\Sigma}^{-1}+\lambda^{-1}(\Theta-\Gamma)\succ 0 \rbrace .
\end{align*}

\begin{lemm}\label{lem41}
Let $(\lambda_k, \Gamma_k, \Theta_k )_{k\in \mathbb{N}}$ be a sequence of elements in $\mathcal{C}_2$ such that
$$
\lim_{k\rightarrow \infty}
\lambda_k = 0.$$
Then $(\lambda_k, \Gamma_k, \Theta_k )_{k\in \mathbb{N}}$ is not an infimizing sequence for $\tilde{J}$ .
\end{lemm}
\proof
We consider two cases separately. 
Let us first analyze the case of sequences $( \lambda_k, \Gamma_k, \Theta_k )$ in which, beside $\lambda_k\rightarrow 0$, we also have $\Vert \lambda_k^{-1}(\Theta_k -\Gamma_k)\Vert\rightarrow \infty$ as $k\rightarrow \infty$.
This implies that the largest singular value of $\lambda_k^{-1}(\Theta_k -\Gamma_k)$ tends to infinity and this, by symmetry, implies in turn that 
\beq\label{eigdivergent}
\lim_{k\rightarrow \infty} \max_{\alpha_k\in \sigma(\lambda_k^{-1}(\Theta_k -\Gamma_k))}
|\alpha_k|=+\infty.
\eeq
We now show that this implies 
\beq\label{meigtominf}
\lim_{k\rightarrow \infty} \min_{\alpha_k\in \sigma(\lambda_k^{-1}(\Theta_k -\Gamma_k))}
\alpha_k=-\infty.
\eeq
To this end, we observe that from (\ref{eigdivergent}) it follows that at least one of the following statements is true:

(\ref{meigtominf}) holds (and in this case we are done) or 
\beq\label{meigtopinf}
\lim_{k\rightarrow \infty} \max_{\alpha_k\in \sigma(\lambda_k^{-1}(\Theta_k -\Gamma_k))}
\alpha_k=+\infty.
\eeq
In the latter case, we use the fact that $\Gamma_k\succeq 0$ and $\lambda_k>0$, so that
\beq\label{maggeig}
\max_{\alpha_k\in \sigma(\lambda_k^{-1}\Theta_k)}\alpha_k
\geq 
\max_{\alpha_k\in \sigma(\lambda_k^{-1}(\Theta_k -\Gamma_k))}
\alpha_k
\eeq
which, together with (\ref{meigtopinf}) gives
\beq\label{meigtominf-mod}
\lim_{k\rightarrow \infty} \max_{\alpha_k\in \sigma(\lambda_k^{-1}\Theta_k)}
\alpha_k=+\infty.
\eeq
Since $\tr(\lambda_k^{-1}\Theta_k)=0$, (\ref{meigtominf-mod}) implies that
\beq\label{mmeigtominf}
\lim_{k\rightarrow \infty} \min_{\alpha_k\in \sigma(\lambda_k^{-1}\Theta_k)}
\alpha_k=-\infty.
\eeq
Now we use again the fact that $\Gamma_k\succeq 0$ and $\lambda_k>0$, so that
\beq
\min_{\alpha_k\in \sigma(\lambda_k^{-1}(\Theta_k -\Gamma_k))} \alpha_k \leq
\min_{\alpha_k\in \sigma(\lambda_k^{-1}\Theta_k)}
\alpha_k
\eeq
which, together with (\ref{mmeigtominf}) implies (\ref{meigtominf}).
In conclusion, for sequences $( \lambda_k, \Gamma_k, \Theta_k )$ of this type and for a sufficiently large $k$,  $\hat{\Sigma}^{-1}+\lambda_k^{-1}(\Theta_k-\Gamma_k)$ is no longer positive definite and therefore these sequences  does not belong to $\mathcal{C}_2$.\\

Second, we consider the case of sequences $( \lambda_k, \Gamma_k, \Theta_k )$ in which, beside $\lambda_k\rightarrow 0$, we also have $\Vert \lambda_k^{-1}(\Theta_k -\Gamma_k)\Vert\rightarrow c$ as $k\rightarrow \infty$, where $c<+\infty$ is a non-negative value.

In this case, it is {\agu not difficult to see that $\forall \varepsilon>0$, $\exists$ $\bar{k}$ such that 
the dual functional satisfies  $\tilde{J}(\lambda_k,\Gamma_k,\Theta_k)>-\varepsilon$, $\forall k\geq \bar{k}$.
In fact, since $\Vert \lambda_k^{-1}(\Theta_k -\Gamma_k)\Vert$ is bounded, there exists $l_0>0$ such that
$\lambda_k^{-1}(\Theta_k -\Gamma_k)\leq l_0I$ for all $k$. Therefore, there exists $l_1>0$ such that
for all $k$, $\hat{\Sigma}^{-1} + \lambda_k^{-1}(\Theta_k -\Gamma_k)\leq l_1I$ and hence there exists $l_2>0$
such that
for all $k$, $\Big{|}\hat{\Sigma}^{-1} + \lambda_k^{-1}(\Theta_k -\Gamma_k)\Big{|}\leq l_2$.
In turn, there exists $l_3\in\mathbb{R}$
such that
for all $k$, $\log\Big{|}\hat{\Sigma}^{-1} + \lambda_k^{-1}(\Theta_k -\Gamma_k)\Big{|}\leq l_3$ and
$-\log\Big{|}\hat{\Sigma}^{-1} + \lambda_k^{-1}(\Theta_k -\Gamma_k)\Big{|}\geq -l_3$.
Eventually, there exists a real constant $l_4:=-l_3-\log\Big{|}\hat{\Sigma}\Big{|}+\delta$
such that, for all $k$, $\tilde{J}(\lambda_k,\Gamma_k,\Theta_k)\geq \lambda_k l_4$.
Since $l_4$ is constant, the the right-hand side of this inequality converges to zero so that, by definition 
$\forall \varepsilon>0$, $\exists$ $\bar{k}$ such that $\lambda_k l_4>-\varepsilon$ $\forall k\geq \bar{k}$.
As a consequence,  $\tilde{J}(\lambda_k,\Gamma_k,\Theta_k)>-\varepsilon$, $\forall k\geq \bar{k}$.}
It is therefore sufficient to exhibit a triple $(\bar{\lambda},\bar{\Gamma},\bar{\Theta})\in \mathcal{C}_2$ for which the dual functional is negative to conclude that sequences 
$( \lambda_k, \Gamma_k, \Theta_k )$ of this kind cannot be minimizing sequences.
Let us consider $(\bar{\lambda},\bar{\Gamma},\bar{\Theta})$ such that $\bar{\lambda}>0$, $\bar{\Gamma}= 0$ and  
\[\bar{\Theta}= -\bar{\lambda}\ofd(\hat{\Sigma}^{-1}).\]
For $\bar \lambda$ sufficiently large, but finite, it is immediate to check that this triple is in $\mathcal{C}_2$.
For this choice of the multipliers and taking into account   \eqref{delta_MAX} we have that
\begin{align*}
\tilde{J}(\bar{\lambda},\bar{\Gamma},\bar{\Theta}) = & -\bar{\lambda}\log\Big{|}\hat{\Sigma}^{-1}+\bar{\lambda}^{-1}(\bar{\Theta}-\bar{\Gamma})\Big{|}
-\bar{\lambda}\log\Big{|}\hat{\Sigma}\Big{|}+\bar{\lambda}\delta\\
 = &-\bar{\lambda}\log\Big{|}(\hat{\Sigma}^{-1}+\bar{\lambda}^{-1}\bar{\Theta})\hat{\Sigma}\Big{|}+\bar{\lambda}\delta\\
= & -\bar{\lambda}\log| \dd(\hat{\Sigma}^{-1})\hat{\Sigma}|+\bar{\lambda}\delta\\
= &-\bar{\lambda}(\delta_{max}-\delta)<0.
\end{align*}

This is sufficient to conclude the proof. In fact, the only other possible case is the one in which
$\lim_{k\rightarrow \infty} \Vert \lambda_k^{-1}(\Theta_k -\Gamma_k)\Vert$ does not exist.
In this case however, we can consider a sub-sequence $( \lambda_{k_j}, \Gamma_{k_j}, \Theta_{k_j} )$ for which the corresponding limit does exist (finite or infinite) and we are thus reduced to one of the previous two cases.
\qed

As a consequence of the previous result we have that the minimization of the dual functional over the set $\mathcal{C}_2$ is equivalent to minimization over the  set:
\begin{align*}
\mathcal{C}_{3}:= &\lbrace (\lambda,\Gamma,\Theta):\lambda\geq \varepsilon , I+\Gamma-\Theta\succeq 0, \Gamma\succeq 0, \Gamma\in\mathbf{D}_n,\\
& \Theta \in \mathbf{M}_n , \hat{\Sigma}^{-1}+\lambda^{-1}(\Theta-\Gamma)\succ 0 \rbrace &
\end{align*}
for a certain $\varepsilon>0$.\\

\noindent The next result provides an upper bound for $\lambda$.

\begin{lemm}
Let $(\lambda_k, \Gamma_k, \Theta_k )_{k\in \mathbb{N}}$ be a sequence of elements in $\mathcal{C}_3$ such that
\beq
\label{ltinf}
\lim_{k\rightarrow \infty}
\lambda_k= \infty.
\eeq
Then $(\lambda_k, \Gamma_k, \Theta_k )_{k\in \mathbb{N}}$ is not an infimizing sequence for $\tilde{J}$. 
\end{lemm}
\proof
Let us consider a sequence $( \lambda_k, \Gamma_k,\Theta_k )_{k\in \mathbb{N}}$ such that (\ref{ltinf})
holds.
 
It follows from the condition $\Theta_k-\Gamma_k \preceq I$ that 
\[\lambda_k^{-1}(\Theta_k-\Gamma_k) \preceq \lambda_k^{-1}I\]
which implies that
\begin{equation}
\begin{aligned}
\tilde{J}(\lambda_k, \Gamma_k, \Theta_k)& = \lambda_k(\log|(\hat{\Sigma}^{-1}+\lambda_k^{-1}(\Theta_k-\Gamma_k))^{-1}\hat{\Sigma}^{-1}|+\delta) \\
& \geq \lambda_k(\log|((\hat{\Sigma}^{-1}+\lambda_k^{-1}I)^{-1}\hat{\Sigma}^{-1})|+\delta)\\ 
& \longrightarrow +\infty
\end{aligned}
\end{equation}
so that $(\lambda_k, \Gamma_k, \Theta_k )_{k\in \mathbb{N}}$ cannot be an infimizing  sequence.
\qed


As a consequence of the previous result, the set $\mathcal{C}_3$ can be further restricted to the set:
\begin{align*}
\mathcal{C}_{4}:=&  \lbrace (\lambda, \Gamma, \Theta):\varepsilon \leq \lambda\leq M , I+\Gamma-\Theta\succeq 0, \Gamma\succeq 0,\\ 
& \Gamma\in\mathbf{D}_n, \Theta \in \mathbf{M}_n , \hat{\Sigma}^{-1}+\lambda^{-1}(\Theta-\Gamma)\succ 0 \rbrace 
\end{align*}
for a certain $M<\infty$.\\

\noindent The next result provides an upper bound for $\Theta-\Gamma$.


\begin{lemm}
Let $(\lambda_k, \Gamma_k, \Theta_k )_{k\in \mathbb{N}}$ be a sequence of elements in $\mathcal{C}_4$ such that
\beq
\label{tmgtinf}
\lim_{k\rightarrow \infty}
\Vert \Theta_k-\Gamma_k\Vert = +\infty.
\eeq
Then $(\lambda_k, \Gamma_k, \Theta_k )_{k\in \mathbb{N}}$ is not an infimizing sequence for $\tilde{J}$. 
\end{lemm}

\proof
From (\ref{tmgtinf}) if follows that the largest singular value of $ (\Theta_k-\Gamma_k)$ tends to $+\infty$ as $k\rightarrow\infty$.
This in turn implies that, as $k\rightarrow\infty$, at least one of the eigenvalues of $ (\Theta_k-\Gamma_k)$ diverges, because  $(\Theta_k-\Gamma_k)$ is symmetric so that its singular values are the absolute values of its eigenvalues.
As before, since $(\Theta_k-\Gamma_k)\preceq I$ holds, the diverging eigenvalues have to tend to $-\infty$. This implies that also  $\hat{\Sigma}^{-1}+\lambda_k^{-1}(\Theta_k-\Gamma_k)$ has an eigenvalue which tends to $-\infty$ as $k\rightarrow\infty$. But, this cannot be the case, because we have the positive definiteness constraint on $\hat{\Sigma}^{-1}+\lambda_k^{-1}(\Theta_k-\Gamma_k)$.
\qed

It follows from the previous result that  there exists $\rho$ such that $|\rho|<\infty$ and 
\[ \Theta -\Gamma\succeq \rho I.\]
Therefore, the set $\mathcal{C}_4$ can be further restricted to the set:\\
\begin{align*}
\mathcal{C}_{5}:= & \lbrace (\lambda,\Gamma,\Theta): \varepsilon \leq \lambda\leq M , \rho I\preceq \Theta -\Gamma\preceq I, \Gamma\succeq 0, \\
& \Gamma\in\mathbf{D}_n, \Theta \in \mathbf{M}_n , \hat{\Sigma}^{-1}+\lambda^{-1}(\Theta-\Gamma)\succ 0 \rbrace. 
\end{align*}


Now observe that in $\mathcal{C}_{5}$ $\Theta$ and $\Gamma$ are orthogonal so that 
if $(\lambda_k, \Gamma_k, \Theta_k )_{k\in \mathbb{N}}$ is a sequence of elements in $\mathcal{C}_5$ such that
\beq
\label{gtinf}
\lim_{k\rightarrow \infty}
\Vert \Gamma_k \Vert = +\infty
\eeq
or
\beq
\label{ttinf}
\lim_{k\rightarrow \infty}
\Vert \Theta_k \Vert = +\infty
\eeq
then (\ref{tmgtinf}) holds.
Then we have the following Corollary.

\begin{corr}
Let $(\lambda_k, \Gamma_k, \Theta_k )_{k\in \mathbb{N}}$ be a sequence of elements in $\mathcal{C}_5$ such that (\ref{gtinf}) or (\ref{ttinf}) holds. Then $(\lambda_k, \Gamma_k, \Theta_k )_{k\in \mathbb{N}}$ is not an infimizing sequence for $\tilde{J}$.
\end{corr}

Thus minimizing over the set $\mathcal{C}_{5}$ is equivalent to minimize over:
\begin{align*}
\mathcal{C}_{6}:= & \lbrace (\lambda, \Gamma, \Theta): \varepsilon \leq\lambda\leq M , \rho I\preceq \Theta -\Gamma\preceq I, 0\preceq\Gamma\preceq \alpha I, \\
& \Gamma\in\mathbf{D}_n,\Theta \in \mathbf{M}_n , \hat{\Sigma}^{-1}+\lambda^{-1}(\Theta-\Gamma)\succ 0 \rbrace 
\end{align*}
for a certain  $\alpha$  such that $0<\alpha<+\infty$.\\

Finally, let us consider a sequence $(\lambda_k, \Gamma_k, \Theta_k )_{k\in \mathbb{N}}\in\mathcal{C}_6$ such that, as $k\rightarrow\infty$, the minimum eigenvalue of $\hat \Sigma+\lambda_k^{-1}(\Theta_k-\Gamma_k)$ tends to zero. This implies that  $\big{|}\hat{\Sigma}^{-1}+\lambda_k^{-1}(\Theta_k-\Gamma_k)\big{|}\rightarrow 0$ and hence $\tilde{J}\rightarrow +\infty$. Thus, such sequence does not infimize the dual functional.
Thus, the final set $\mathcal{C} $ is 
\begin{align*} 
\mathcal{C}:= & \lbrace (\lambda, \Gamma, \Theta): \varepsilon\leq \lambda\leq M , \rho I\preceq \Theta -\Gamma\preceq I, 0\preceq\Gamma\preceq \alpha I, \\
& \Gamma\in\mathbf{D}_n,\Theta \in \mathbf{M}_n , \hat{\Sigma}^{-1}+\lambda^{-1}(\Theta-\Gamma)\succeq \beta I \rbrace 
\end{align*}
for a suitable $\beta >0$.

 Summing up we have the following Theorem. 
\begin{teor}
Problem (\ref{Dual}) is equivalent to
\begin{equation}
\label{dual_compact}
\begin{aligned}
\min_{(\lambda,\Gamma, \Theta) \in \mathcal{C}} \tilde{J}(\lambda,\Gamma, \Theta).
\end{aligned}
\end{equation}
Both these problems admit solution.
\end{teor}
\proof
Equivalence of the two problems has already been proven by the previous argument.
Since $\mathcal{C}$ is closed and bounded and, hence, compact, and $\tilde J$ is continuous over $\mathcal{C}$, by the Weierstrass's Theorem the minimum exists.
\qed

Before discussing the uniqueness of the solution to (\ref{Dual}), it is convenient to further simplify the dual optimization problem: consider the function
  
\[F(\lambda,X):=-\lambda [\log(|\hat{\Sigma}^{-1}+\lambda^{-1}X|)
+\log|\hat{\Sigma}|-\delta]
\]
where $\lambda>0 $ and $X\in {\mathbf Q}_n$. Note that 
$$
F(\lambda,\Theta-\Gamma)=\tilde J (\lambda,\Gamma,\Theta).
$$
Moreover, $\Theta$ and $\Gamma$ are orthogonal over $\cal C$
 so that minimizing $\tilde J$ over ${\cal C}_0$ is equivalent to  minimize $F$
over the corresponding set
\begin{align*} 
\mathcal{C}_F:=  \lbrace (\lambda, X): & \  \lambda>0, \ 
X\in \mathbf{Q}_n, X\preceq I,  \\
& \ -\dd(X)\succeq 0,  \ \hat{\Sigma}^{-1}+\lambda^{-1}X \succ 0 \rbrace.
\end{align*}
Therefore, from now on we can consider the following problem
\begin{equation}
\label{dual_simplif}
\begin{aligned}
\min_{(\lambda,X) \in \mathcal{C}_F} F(\lambda,X).
\end{aligned}
\end{equation}
Once obtained the optimal solution $(\lambda^*,X^*)$ we can recover the optimal values of the original multipliers simply by setting
$\Theta^*=\ofd(X^*)$ and $\Gamma^*=-\dd(X^*)$.

\subsection{Uniqueness of the solution of the dual problem}
The aim of this Section is to show that Problem (\ref{dual_simplif}) (and, hence
Problem (\ref{Dual})) admits a unique solution.
Since $\tilde J$ is the opposite of the dual objective function, $\tilde J$ is convex over $\cal C$.
It is then easy to check that $F$ is also a convex function over the convex set
$\mathcal{C}_F$. 
However, as we will see, $F$ is not strictly convex. Accordingly, establishing the uniqueness of the minimum is not a trivial task.

The following Proposition, whose proof is in Appendix, characterizes the second variation of $F$ in  direction
$(\delta \lambda, \delta X)$, i.e. $\delta^{2}F (\lambda,X;\delta \lambda, \delta X)$.
\begin{propo} \label{prop_hessian}
Let  $x:=\text{vec}(X)$, $\delta x:=\text{vec}(\delta X)$, and $K:= (\hat{\Sigma}^{-1}+\lambda^{-1}X)^{-1}\otimes (\hat{\Sigma}^{-1}+\lambda^{-1}X)^{-1}$. Let also
\[ H:= \left[ \begin{array}{cc}
\lambda^{-3}x\tp K x & -\lambda^{-2}x\tp K  \\
-\lambda^{-2}Kx & \lambda^{-1}K\\
\end{array} \right]\in \mathbb{R}^{(1+n^2) \times (1+n^2)}. \]
Then, we have
$$
\delta^{2}F (\lambda,X;\delta \lambda, \delta X)=[\delta \lambda\ \ \delta x\tp]
H\bmat{c}\delta \lambda\\ \delta x\emat.
$$
\end{propo}

Since in $\mathcal{C}_F$ we have that $K	\in {\mathbf Q}_n$ is positive definite and $\lambda>0$, the matrix  $H$,  which has clearly the meaning of the Hessian of $F$, has at least rank equal to $n^2$. Moreover, $Hw=0$ with $w=[\,\lambda \; x\tp \,]\tp $. We conclude that $H$ has rank equal to $n^2$.

 This means that $F$ is convex and there is exactly one direction along which 
$F$ is not strictly convex.
We now analyse this direction in the neighbourhood of the optimal solution.
\begin{lemm}\label{optonbound}
Any optimal solution $(\lambda^{*},X^{*})$ minimizing $F$ over $\mathcal{C}_F$  lies on the boundary of 
$\mathcal{C}_F$ and, specifically, is such that $I-X^*$ is singular.
\end{lemm}
\proof
Let $(\lambda^{*},X^{*})$ be an optimal solution and assume, by contradiction, that
$(\lambda^{*},X^{*})$ does not belong to the boundary of the feasible set $\mathcal{C}_F$, so that, in particular, $X^*\prec I$. Thus there exists $\varepsilon>0$ such that 
$(1+\varepsilon)X^*\prec I$ so that $$((1+\varepsilon)\lambda^{*},(1+\varepsilon)X^{*})\in \mathcal{C}_F.$$
Now a direct computation yields
\beq
F((1+\varepsilon)\lambda^{*},(1+\varepsilon)X^{*})=(1+\varepsilon)F(\lambda^{*},X^{*})<F(\lambda^{*},X^{*})
\eeq
where the last inequality is a consequence of the fact that, as 
we have already seen in the proof of Lemma \ref{lem41}, the optimal value of $\tilde{J}$ (and, hence, of $F$) is negative.
This a contradiction as $F(\lambda^{*},X^{*})$ is assumed to be a minimum. 
\qed

\begin{rem}
Notice that for any $(\lambda_0,X_0)\in \mathcal{C}_F$, the direction $(\varepsilon\lambda_0,\varepsilon X_0) $ (which, by the way, is the direction considerered in Lemma
\ref{optonbound} for the specific case of the optimal solution $(\lambda^{*},X^{*})$)
is exactly the unique direction along which $F$ is not strictly convex. In fact,
along this direction $F$ is clearly a linear function of $\lambda$. 
Notice also that $F$ is constant along this direction if and only if  $F(\lambda_0,X_0)=0$.
Since at any optimal solution $(\lambda^{*},X^{*})$ $F$ is necessarily negative,  $F$ is not constant along the direction $(\varepsilon\lambda^{*},\varepsilon X^{*}) $ (which is the only direction along which $F$ is not strictly convex). 
\end{rem}
As a consequence of this observation, we have the following result.
\begin{corr}\label{corstconvex}
Let $(\lambda_0,X_0)$ be a given point in $\mathcal{C}_F$.
If $w:=(\delta\lambda,\delta X)$ is any direction along which $F(\lambda_0,X_0)$ is constant, i.e.
$F(\lambda_0,X_0)=F(\lambda_0+\alpha \delta\lambda,X_0+\alpha \delta X)$ for any $\alpha$ such that $|\alpha|>0$ is sufficiently small, then $F(\lambda_0,X_0)=0$.
\end{corr}
We are now ready to prove our main result.
\begin{teor}
The dual problem admits a unique solution.
\end{teor}
\proof
Assume, by contradiction, that there are two optimal solutions $(\lambda^{*}_1,X_1^{*})$ and $(\lambda^{*}_2,X_2^{*})$. By the convexity of $\mathcal{C}_F$, the whole segment $S$ connecting $(\lambda^{*}_1,X_1^{*})$ to $(\lambda^{*}_2,X_2^{*})$ belongs also to $\mathcal{C}_F$. 
It follows by the convexity of $F(\cdot, \cdot)$ that all the points in $S$ are optimal solutions. Notice, in passing, that  in view of Lemma \ref{optonbound}, this implies that $S$  belongs to the boundary of $\mathcal{C}_F$.
Now, $F$ is clearly  negative and constant  along $S$ and this is a contradiction in view of  Corollary \ref{corstconvex}. 
\qed

\section{Recovering the solution of the primal problem} \label{sec_recovery}
By the uniqueness of the solution of the dual problem we know that the duality gap between the primal and the dual problem is zero. This allows us to recover the solution of the primal problem.\\
First, the optimal  $\Sigma$ can be easily recovered by substituting the optimal solution of the dual problem $(\lambda^*,\Theta^*,\Gamma^*)$ into \eqref{der_var_wrt_sigma_opt}. Recovering the optimal $L$ is slightly more involved; since the duality gap is zero, from the KKT conditions we have:

\begin{equation}
\label{KKT_1}
\tr(\Lambda L)=0
\end{equation}
\begin{equation}
\label{KKT_2}
\tr(\Gamma(\Sigma-L))=0
\end{equation}
\begin{equation}
\label{KKT_3}
\tr(\Theta(\Sigma-L))=0.
\end{equation}

We begin by considering \eqref{KKT_1}. It follows from \eqref{der_var_wrt_r} that
\[\Lambda=I+\Gamma-\Theta\]
where we now know that $\Lambda$ has deficient rank. Thus, we consider the following reduced singular value decomposition
\begin{equation}
\label{reduced_svd}
\Lambda=U S U\tp 
\end{equation}
with $S\in \mathbf{Q}_{n-r}$ positive definite, i.e. $n-r$ is the rank of $\Lambda$, and $U\in\mathbb{R}^{n\times n-r}$ such that $U\tp U =I_{n-r}$. We plug \eqref{reduced_svd} in \eqref{KKT_1} and get
\begin{equation}
\begin{aligned}
0=\text{tr}[\Lambda L]= \text{tr}[USU\tp  L]
\Rightarrow U\tp LU=0.
\end{aligned}
\end{equation} 
Then, by selecting a matrix $\tilde{U}\in \mathbb{R}^{n\times r}$ whose columns form an orthonormal bases of $[{\rm im}(U)]^\bot$, we can
express  $L$ as:
\begin{equation}
\label{reduced_sdv_R}
L=\tilde{U}Q\tilde{U}\tp 
\end{equation}
with $Q\in \mathbf{Q}_r$. Note that, in view of the fact that the columns of $\tilde{U}$ form the orthogonal complement of the image of $U$,the relationship $U\tp \tilde{U}=0$ holds.\\

By (\ref{KKT_3}), we know that $\Sigma-L$ is diagonal. Thus, we plug \eqref{reduced_sdv_R} into (\ref{KKT_3}) and obtain a linear system of equations: $\ofd(\Sigma-\tilde{U}Q\tilde{U}\tp )=0$, or equivalently, 
\begin{equation}
\begin{aligned}
\ofd(\tilde{U}Q\tilde{U}\tp )=\ofd(\Sigma).
\end{aligned}
\end{equation}
In an analogous fashion, using (\ref{KKT_2}) we obtain an additional system of linear equations.
In virtue of the fact that both the dual and the primal problem admit solution the resulting system of equations always admits solution in $Q$. Moreover, the solution of this system of equations is unique if and only if the solution of the primal problem is unique.

\section{Numerical Implementation}\label{algorithm}

We propose an algorithm for finding the numerical solution of Problem (\ref{dual_simplif}). First, recall that the optimal solution lies in the boundary characterized by constraints $-\dd (X)\succeq 0$ and $X\preceq I$. Finding a descending direction $(\lambda,X)$ for $F(\lambda,X)$ satisfying simultaneously these two constraints is not trivial. Then we resort to the Alternating Direction Method of Multipliers (ADMM) algorithm, \cite{boyd2011distributed}, for decoupling such constraints. Then, the corresponding ADMM updates can be performed by using a projection gradient algorithm. To this end, we rewrite (\ref{dual_simplif}) by introducing the new variable $Y\in \mathbf{Q}_n$ defined as $Y:=I-X$ :
\begin{equation*}
\begin{aligned}
\underset{\lambda,X,Y}{\min} &F(\lambda,X) \\
\hbox{subject to }& {(\lambda,X)\in \mathcal{C}^{*}_{\lambda, X} },\;\; {Y\in \mathcal{C}^{*}_Y}\\
& Y=I-X
\end{aligned}
\end{equation*}
with $\mathcal{C}^{*}_{\lambda, X}$, and $\mathcal{C}^{*}_{Y}$ defined, respectively, as

\begin{alignat*}{3}
\mathcal{C}^{*}_{\lambda, X} := & \lbrace (\lambda, X): \lambda>0, X\in\mathbf{Q}_n,  \hat{\Sigma}^{-1}+\lambda^{-1}X \succ 0, \\
& \quad  -\dd(X)\succeq 0 \rbrace \\
\mathcal{C}^{*}_{Y} := & \lbrace Y: Y\in \mathbf{Q}_n, Y\succeq 0 \rbrace .
\end{alignat*}
The {\em augmented Lagrangian} (see \cite{boyd2011distributed}) for the problem is
\begin{equation*}
\mathcal{L}_{\rho} (\lambda, X, Y, M) =F(\lambda,X)+\langle M, Y-I+X \rangle +\dfrac{\rho}{2}\Vert Y-I+X \Vert^{2}_F
\end{equation*}
where $M\in \mathbf{Q}_n$. Accordingly, given the initial values $\lambda^0$, $X^0$, $Y^0$ and $M^0$, the ADMM updates are:
\begin{align}
\begin{split}\label{1admmupdate}
(\lambda^{(k+1)},X^{(k+1)}) {}& :=\argmin_{(\lambda, X) \in \mathcal{C}^{*}_{\lambda, X}} \mathcal{L}_{\rho} (\lambda, X, Y^{(k)}, M^{(k)})
\end{split}\\
\begin{split}\label{2admmupdate}
Y^{(k+1)} {}& :=\argmin_{Y \in \mathcal{C}^{*}_{Y}} \mathcal{L}_{\rho} (\lambda^{(k+1)}, X^{(k+1)}, Y, M^{(k)})
\end{split}\\
\begin{split}\label{3admmupdate}
M^{(k+1)} {}& :=M^{(k)}+\rho(Y^{(k+1)}-I+X^{(k+1)})
\end{split}
\end{align}
where $\rho>0$ is the penalty parameter. Here, we choose $\rho=0.5$.

Problem \eqref{1admmupdate} has not a closed form solution. Thus, the solution is approximated by a projective gradient step:
\begin{align*}
\lambda^{(k+1)}:= \lambda^{(k)}-t_k \nabla_\lambda \mathcal{L}_\rho (\lambda^{(k)}, X^{(k)}, Y^{(k)}, M^{(k)})\\
X^{(k+1)}:= \Pi_{\mathcal{C}^{*}_{X}}(X^{(k)}-t_k \nabla_X \mathcal{L}_\rho(\lambda^{(k)}, X^{(k)}, Y^{(k)}, M^{(k)})
\end{align*}
where $\Pi_{\mathcal{C}^{*}_{X}}$ denotes the projector operator from $\mathbf{Q}_n$ onto 
\alg{\mathcal{C}^{*}_{X}:=\lbrace X: X\in\mathbf{Q}_n, -\dd(X)\succeq 0 \rbrace,\nn} and $\nabla_X \mathcal{L}_\rho$, $\nabla_X \mathcal{L}_\rho$ denotes the gradient with respect to $\lambda$ and $X$, respectively:

\begin{equation*}
\begin{aligned}
\nabla_X \mathcal{L}_\rho(\lambda, X,& Y, M):=\\ 
& -\log |\hat{\Sigma}^{-1}+\lambda^{-1} X | -\log|\hat{\Sigma}|+ \delta \\
&+\lambda^{-1}\tr\big((\hat{\Sigma}^{-1}+\lambda^{-1}X)^{-1}X\big) \\
\nabla_X \mathcal{L}_\rho (\lambda, X,& Y, M):= \\ 
&-(\hat{\Sigma}^{-1}+\lambda^{-1}X)^{-1}+M+\rho(Y-I-X).
\end{aligned}
\end{equation*}
It is not difficult to see that
\alg{[\Pi_{\mathcal{C}^{*}_{X}}(A)]_{ij}=\left\{\begin{array}{ll}0, & \hbox{if $i=j$ and $[A]_{ij}>0$} \\ \hbox{$[A]_{ij}$}, & \hbox{otherwise} \end{array}\right.\nn}
where $[A]_{ij}$ denotes the entry in position $(i,j)$ of matrix $A\in\mathbf{Q}_n$.
The step size $t_k$ is determined at each step $k$ in a iterative fashion: we start by setting $t_k=1$ and we decrease it progressively until the two conditions $\lambda^{(k+1)}>0$ and $\hat{\Sigma}^{-1}+\lambda^{-1}X \succ 0$ are met and the so-called Armijo condition, \cite{BOYD_CONVEX_OPTIMIZATION}, are satisfied.\\

Problem \eqref{2admmupdate} can be rewritten as 
\[Y^{(k+1)}=\argmin_{Y\in\mathcal{C}^{*}_{Y}} \Vert I-X^{(k+1)}-\dfrac{1}{\rho}M^{(k)} -Y \Vert_F.\]
We introduce the projection operator $\Pi_{\mathcal{C}^{*}_{Y}}:\mathbf{Q}_n\rightarrow \mathcal{C}^{*}_{Y}$ 
which is defined as 
\[\Pi_{\mathcal{C}^{*}_{Y}}(W):=\argmin_{Z\in \mathcal{C}^{*}_{Y} }\Vert W-Z \Vert_F^{2}.\]
It is not difficult to see that, if $A=UDU\tp$ is the eigenvalue decomposition of the matrix $A\in\mathbf{Q}_n$, then
\[\Pi_{\mathcal{C}^{*}_{Y}}(A)=U\diag(f(d_1),...,f(d_n))U\tp\]
where
\[f(d_i):= \begin{cases}
		d_i, & \text{if  } d_i\geq 0 \\    0 & \text{otherwise.}
		\end{cases}\]
Then the solution of \eqref{2admmupdate} becomes
\[Y^{(k+1)}=\Pi_{\mathcal{C}^{*}_{Y}}( I-X^{(k+1)}-\dfrac{1}{\rho}M^{(k)}).\]

In order to set the stopping criteria for the algorithm we define the primal and dual residual matrices:
\begin{align*}
R^{(k+1)}& := Y^{(k+1)}-I+X^{(k+1)},\\
S^{(k+1)}& := \rho_k(Y^{(k+1)}-Y^{(k)}).
\end{align*}
The algorithm reaches an acceptable solution when the following conditions are met, \cite{boyd2011distributed}:
\begin{align*}
\Vert R^{(k+1)} \Vert_F & \leq n\epsilon^{abs}+\epsilon^{rel}\max\lbrace\sqrt{n},\Vert X^{k} \Vert_F, \Vert Y^{k} \Vert_F\rbrace,\\
\Vert S^{(k+1)}\Vert & \leq n\epsilon^{abs}+\epsilon^{rel}\Vert M^{(k)}\Vert_F
\end{align*}
where $\epsilon^{rel}=10^{-4}$ and $\epsilon^{abs}=10^{-4}$ are the relative and the absolute tolerance, respectively.

\section{Numerical Examples} \label{simulations}
\subsection{Synthetic data}
In this Section we consider Monte Carlo studies composed by 200 experiments whose structure is as follows. For each experiment:
\begin{itemize}
\item  we consider a factor model having the structure of  (\ref{fact_model}) with the cross sectional dimension is $n=40$; $L$ and $D$ are randomly generated in such a way that $L$ has rank equal to $r$ ({\em a priori} fixed), $D$ is diagonal and the signal-to-noise ratio {\vale (defined as $\Vert L\Vert/\Vert D \Vert$)} between the latent and the idiosyncratic  components is equal to one;
\item a data sequence of length $N$ for $y$ is generated;
\item we compute the sample covariance matrix $\hat \Sigma$ from this data;
\item we compute the estimate $\delta_\alpha$ of $\delta$ using the empirical procedure of Section \ref{sec_choice_delta} with $\alpha=0.5$; 
{\vale
\item we compute the solution $(L_{\scriptscriptstyle OPT},\Sigma_{\scriptscriptstyle OPT})$ of Problem (\ref{primal}) where we replace $\delta$ with $\delta_{\alpha}$. Let $\lambda_i$, $i=1\ldots n$, denote the singular values of $L_{\scriptscriptstyle OPT}$ and define $i_{max}$ as the first $i$ such that $\lambda_{i+1}/\lambda_{1}<0.05$.  Then, we define the ``numerical rank'' of $L_{\scriptscriptstyle OPT}$ as:
\alg{r_{\scriptscriptstyle OPT}:= \max_{i\leq i_{max}} \lambda_{i}/\lambda_{i+1}}
\item we compute the solution of the standard problem (\ref{original_problem}) (exact decomposition with trace heuristic)  and, with the same procedure of the previous point, we compute the numerical rank, $ r_{{\scriptscriptstyle ED}}$, of the corresponding low rank matrix;
\item we compute the minimum number of factors from the data sequence for $y$ by applying the three methods proposed by Bai and Ng \cite{bai2002determining}, namely: ICP1, ICP2 and ICP3. We denote the corresponding estimates by $r_{{\scriptscriptstyle ICP1}}$, $r_{{\scriptscriptstyle ICP2}}$ and $r_{{\scriptscriptstyle ICP3}}$, respectively;
\item we compute the minimum number of factors from the data sequence for $y$ by applying the method proposed by Lam and Yao \cite{lam2012factor}.\footnote{\agu The estimation procedure for this method requires to set a parameter $k_0$ for the selection  which only general considerations are provided: we decided to select $k_0$ using an ``oracle'' procedure i.e. for each Monte Carlo run we choose the value of $k_0$ which yields the most favourable result.} We denote by $r_{{\scriptscriptstyle LY}}$ the resulting estimate of the rank.
}
\end{itemize}
{\vale 
Finally, we compute the root mean squared error:
{\agu\begin{equation}
\label{error}
e=\sqrt{\dfrac{1}{200}\sum_{i=1}^{200}(r_{i}^{\sharp}-r)^2}
\end{equation}
}  
for $r^{\sharp}=\lbrace r_{\scriptscriptstyle OPT}, r_{{\scriptscriptstyle ED}}, r_{{\scriptscriptstyle ICP1}}, r_{{\scriptscriptstyle ICP2}}, r_{{\scriptscriptstyle ICP3}}, r_{{\scriptscriptstyle LY}}\rbrace$ and {\agu where} $r$ is the true rank of the data generating process.
}
{\vale
Table \ref{r4tab}
\begin{table}
\vale \scriptsize
\centering
\begin{tabular}{|l|c|c|c|c|c|c|}
\hline
\rule{0pt}{1.0\normalbaselineskip}
\multirow{2}{*}{ } & \multirow{2}{*}{\shortstack[c]{Proposed\\ method}}  & \multirow{2}{*}{\shortstack[c]{Exact\\ Decomposition} } & 
\multicolumn{3}{c|}{Bai $\&$ Ng} & \multirow{2}{*}{\shortstack[c]{Lam $\&$\\ Yao}}\\
\cline{4-6}
 \rule{0pt}{1.0\normalbaselineskip}
& & & ICP1 & ICP2 & ICP3 & \\
\hline
\rule{0pt}{1.0\normalbaselineskip}

 N=200  & 0.500 & 3.752 & 3.589 & 2.546 & 7.506 & 5.529\\
 
 N=500  & 0.000 & 0.9618 & 2.271 & 2.273 & 4.236 & 5.347\\

 N=1000  & 0.000 & 0.6557 & 3.587 & 3.213 & 3.927 & 5.421\\
\hline
\end{tabular}
\medskip
\caption{Average root mean squared error between the estimated numerical rank and the true rank $r=4$.} 
\label{r4tab}
\end{table} shows error \eqref{error} for three Monte Carlo studies where $r=4$ and the sample size is $N=200$, $N=500$ and $N=1000$, respectively.

Usually, the problem becomes more challenging when the rank $r$ of the data generating process increases (yet remaining below the Ledermann bound).
For this reason we repeat the above Monte Carlo studies for the case $r=10$ (considering again the three sample sizes $N=200$, $N=500$, $N=1000$). The corresponding root mean squared errors \eqref{error} are reported in Table \ref{r10tab}.

 As one can see, {\agu in all these six Monte Carlo studies} the proposed method outperforms the others.}  
\begin{table}
\vale \scriptsize
\centering
\begin{tabular}{|l|c|c|c|c|c|c|}
\hline
\rule{0pt}{1.0\normalbaselineskip}
\multirow{2}{*}{ } & \multirow{2}{*}{\shortstack[c]{Proposed\\ method}}  & \multirow{2}{*}{\shortstack[c]{Exact\\ Decomposition} } & 
\multicolumn{3}{c|}{Bai $\&$ Ng} & \multirow{2}{*}{\shortstack[c]{Lam $\&$\\ Yao}}\\
\cline{4-6}
\rule{0pt}{1.0\normalbaselineskip}
& & & ICP1 & ICP2 & ICP3 & \\
\hline 
\rule{0pt}{1.0\normalbaselineskip}

N=200  & 2.170 & 7.218 & 5.888 & 5.629 & 8.254 & 6.943\\
 
 N=500  & 0.174 & 5.221 & 5.214 & 5.258 & 5.812 & 6.536\\

 N=1000  & 0 & 2.961 & 5.302 & 5.196 & 5.490 & 6.669\\
\hline
\end{tabular}
\medskip
\caption{Average root mean squared error between the estimated numerical rank and the true rank $r=10$.} 
\label{r10tab}
\end{table} 

\medskip

{\vale We now analyze how well the proposed method recovers the subspace of $L$ by considering the following measure of discrepancy. Let $L=AA^{\tp}$ be the low rank matrix of the data generating process and consider the singular value decomposition of $L_{\scriptscriptstyle OPT}$, that is $L_{\scriptscriptstyle OPT}=USV^{\tp}$. 
Let $\tilde{U}:=U_{[1:n,1:r_{\scriptscriptstyle OPT}]}, \; \tilde{U}\in \mathbb{R}^{n\times r_{\scriptscriptstyle OPT}}$ be the matrix formed by the first $r_{\scriptscriptstyle OPT}$ columns of $U$ and $\tilde{S}:= S_{[1:r_{\scriptscriptstyle OPT},1:r_{\scriptscriptstyle OPT}]}, \; \tilde{S}\in\mathbb{R}^{r_{\scriptscriptstyle OPT}\times r_{\scriptscriptstyle OPT}}$ be the top left $r_{\scriptscriptstyle OPT}\times r_{\scriptscriptstyle OPT}$ sub-matrix of $S$.
We define the projector onto the subspace of $A_{\scriptscriptstyle OPT}:=\tilde{U}\tilde{S}$ as
\[P:=A_{\scriptscriptstyle OPT}(A_{\scriptscriptstyle OPT}^{\tp}A_{\scriptscriptstyle OPT})^{-1}A_{\scriptscriptstyle OPT}^{\tp}.\]
Then, a measure of discrepancy between the subspace of $A$ and the subspace of $A_{\scriptscriptstyle OPT}$ is given by:
\begin{equation}
\label{subspace_rec}
s(A_{\scriptscriptstyle OPT}):=\tr(A^{\tp}P A)/ \tr(A^{\tp}A)
\end{equation}
where $s(A_{\scriptscriptstyle OPT})$ takes value between 0 and 1. Note that, if $s(A_{\scriptscriptstyle OPT})=1$ then $A_{\scriptscriptstyle OPT}$ recovers exactly the image of $A$. Figure \ref{fig:subspace} (left hand side panel) shows the box-plots for error \eqref{subspace_rec} in the three Monte Carlo studies. \\
\begin{figure}[htb]
	\centering
	\includegraphics[trim={4.5cm 0 4.5cm 0},clip,width=0.48\textwidth]{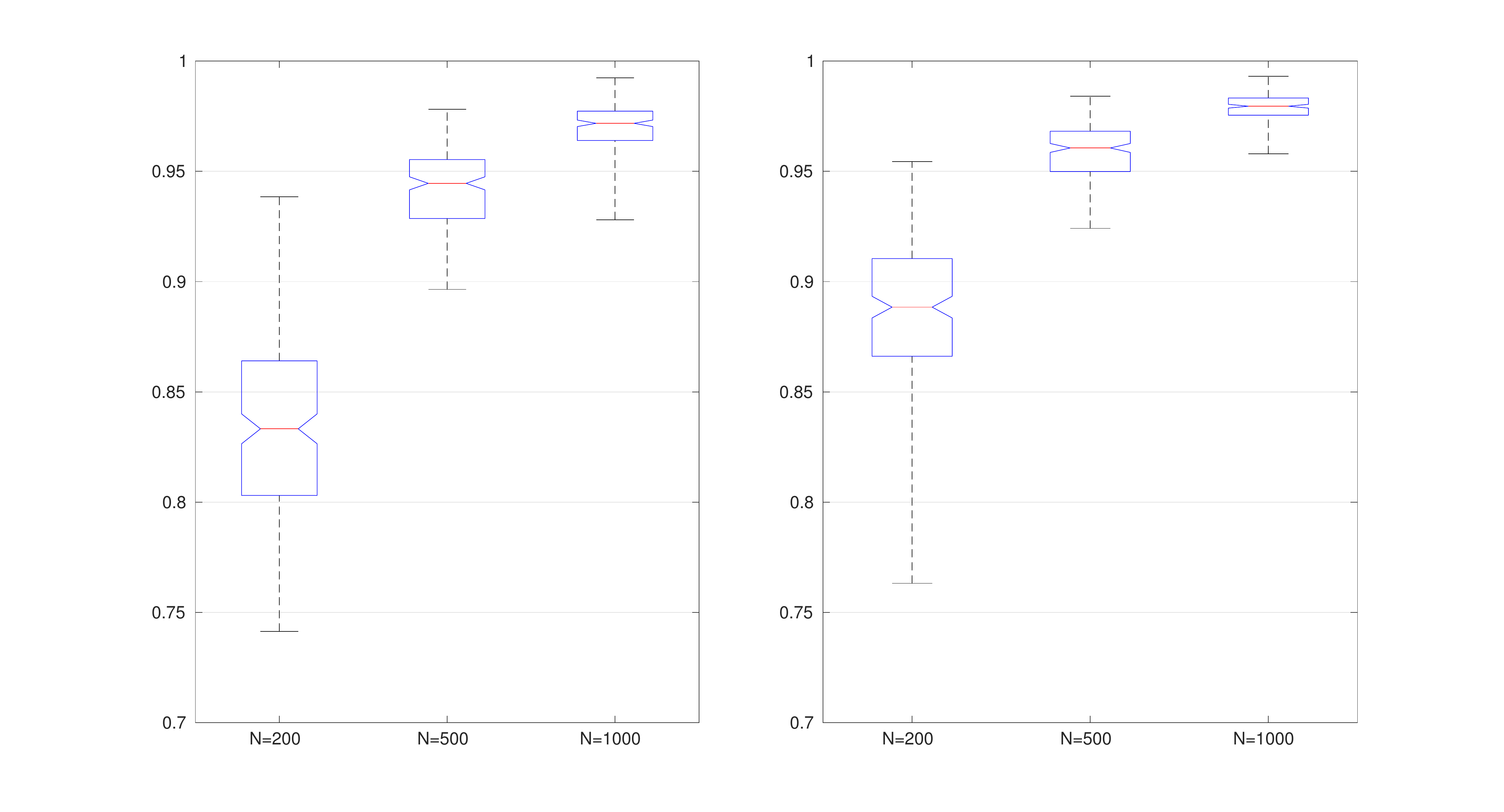}
\caption{\mz Performances of the subspace recovery: box-plots of the quantity \eqref{subspace_rec} for the Monte Carlo studies with sample sizes $N=200, 500, 1000$ and for the case $r=4$ (left hand side panel) and $r=10$ (right hand side panel).}
\label{fig:subspace}
\end{figure}

{\mz Finally, we consider the example illustrated in Figure \ref{fig:barplot1} of the Introduction.
By applying our method, we obtain the situation illustrated in Figure \ref{fig:barplot2} showing  that our approach provides a numerical rank equal to the true value of  $r$.}

\begin{figure}[htb]
	\centering
	\includegraphics[trim={4.5cm 2cm 4.5cm 1.5cm},clip,width=0.46\textwidth]{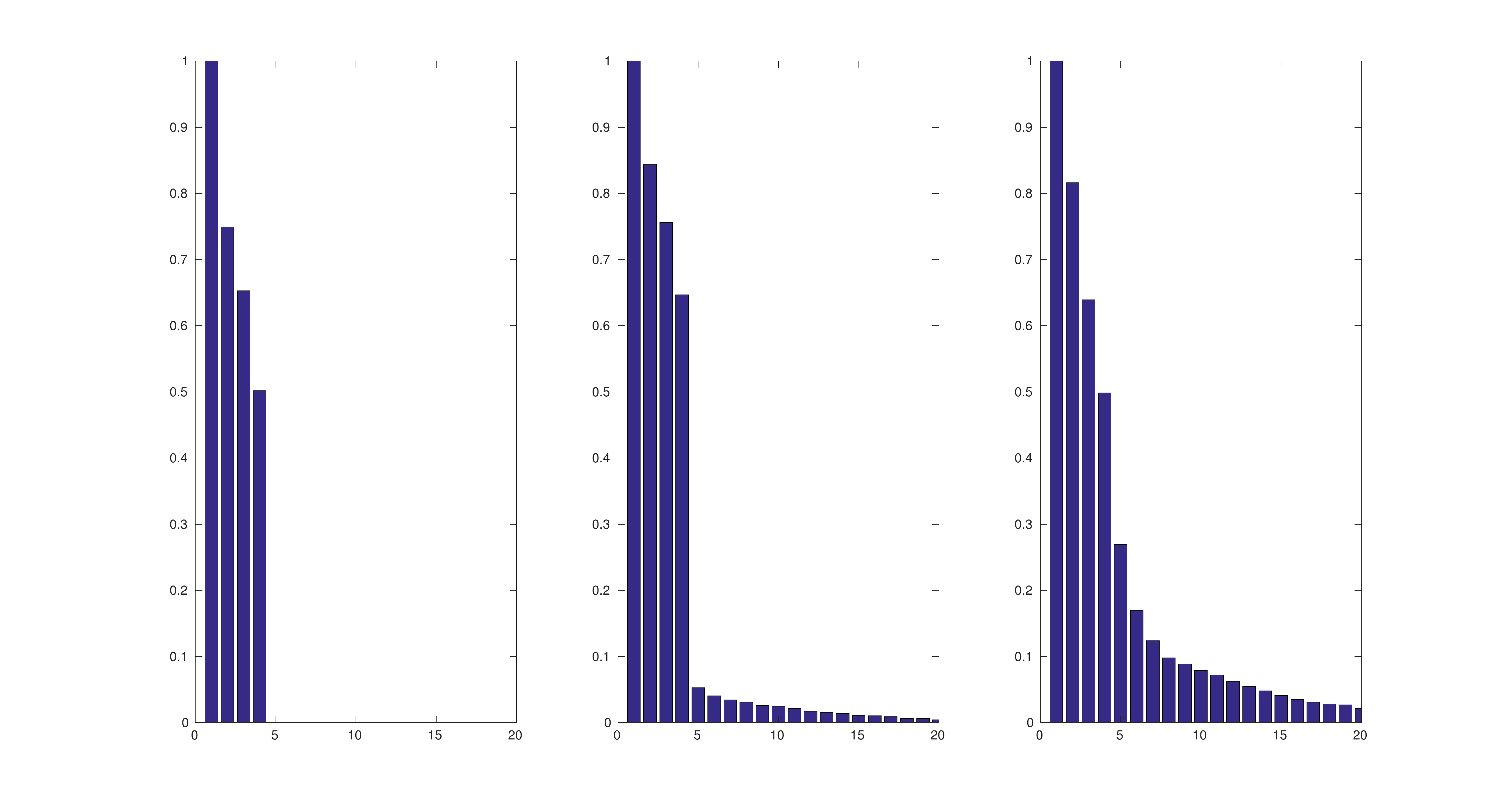}
\caption{{\vale A sample of numerosity 1000 has been generated from a factor model with $n=40$ and $r=4$. The Figure displays the first twenty singular values of the true matrix $L$ (on the left) and of the matrices $L_{\scriptscriptstyle OPT}$ (in the middle) and $L_{\scriptscriptstyle ED}$ (on the right) estimated, respectively, with the proposed method and with the trace heuristic with exact decomposition.}}
\label{fig:barplot2}
\end{figure}

\subsection{Data analysis for investment decision} \label{real_data}

In this sub-section we consider a cross section of $9$ financial indicators ($n=9$) collected across $94$ different sectors ($N=94$) of the US economy (each data vector represents the average for that sector). The data are taken from 
\url{http://www.stern.nyu.edu/~adamodar/pc/datasets/betas.xls} (data downloaded on June 2017). \\ 
The considered financial indicators can be computed from the balance sheet of the companies and from stocks market development and are customarily used for investment decisions. These indicators are: the beta, that is the systemic risk arising from the exposure to general market movements, the debt/equity ratio, the tax rate, the unlevered beta, the cash/firm value ratio, the unlevered beta corrected for cash, the Hi-Lo risk, the standard deviation of equity and the standard deviation of operating income. For a proper description of these indicators we refer to the aforementioned web site while for a more general treatment see for example \cite{koller2010valuation}.\\
It is reasonable to expect that the variability of the listed indicators may be successfully explained by a smaller number of factors and motivated by this {\vale reason} we estimate the sample covariance matrix and we apply the proposed approach (with $\alpha=0.5$).
Indeed, we obtain an estimate $\hat{r}=3$ for the number of latent factors.
These seems to be reasonable, since the common variability of these indicators may be explained by factors such as the general market trend, the different fiscal regime and the different optimal capital structure across sectors.
{\gus In this case, the POET method proposed in \cite{fan2013large} provides an estimate of $1$ latent factor and the methods 
proposed in \cite{bai2002determining} and \cite{lam2012factor} provide all an estimate of $8$ latent factors: the latter number does  not seem very reliable as it is larger  than the upper bound of  $7$ latent factors provided by the method based on exact decomposition of the covariance matrix.}\footnote{\agu The eighth and ninth eigenvalues of the corresponding $L$ matrix obtained with this method are only numerically non-zero as they are smaller than $10^{-14}\lambda_1$, with $\lambda_1$ being the largest eigenvalue of $L$.}

\section{Conclusion and further research directions}\label{sec_concl}
In this paper we have proposed a {\mz new method to estimate the number of factors} for the realistic situation in which the
covariance matrix of the data is estimated with an error that is not negligible. 
{\vale 

A question which arises naturally concerns the statistical properties of the proposed estimator, and, in particular, its asymptotic properties as the sample size approaches infinity. 
{\gus This is a complex issue that certainly cannot be fully addressed in the context of the present paper.
We present only some ideas in this direction and a heuristic road-map that can be followed.}
A primary issue is the rank consistency of the minimum trace estimator, note that a similar matter has been studied e.g. in \cite{bach2008consistency} for the Lasso problem. 
Restricting to the cases in which the minimizers of the minimum rank problem and of the minimum trace problem coincide, a
{\gus possible} argument may be the following. By the construction presented in Section \ref{sec_choice_delta}, it holds that
\[Pr[2 \mathcal{D}_{KL}(\Sigma \Vert \hat{\Sigma})< \delta]= \alpha.\]
{\gus We can} let the desired precision $\alpha$ be a function of the sample size $N$, {\gus and choose $\alpha(N)$ such that, as $N \rightarrow \infty$, it holds that $\alpha(N)\rightarrow 1 $. Moreover, 
 we let $\alpha(N)\rightarrow 1 $ sufficiently slowly so that it is reasonable to expect that } $\delta(\alpha(N))\rightarrow 0$ {\mz because $\hat\Sigma \rightarrow \Sigma$ almost surely}.
 {\gus Consequently, as $N \rightarrow \infty$, the neighbourhood of $\hat{\Sigma}$ in which we seek for the solution 
becomes smaller and smaller and it contains the ``true'' $\Sigma$ with probability tending to $1$.
Moreover, the minimum rank problem \eqref{rank_min}  is} a lower semi-continuous function of $\Sigma$ (see \cite{ning2015linear}, Proposition 1) and, being integer valued, it does not decrease in a sufficiently small neighborhood of $\Sigma$. Therefore, it seems reasonable to conclude that $\forall \, \epsilon > 0 \; \exists \; \bar{N}$ such that $\forall N> \bar{N}$ it holds that
\[Pr[r_{\scriptscriptstyle OPT} \neq r_{true}]< \epsilon.\]
A rigorous study of this heuristic argument will be subject of future investigation.

Another natural direction of research, which will be subject of future investigation, is the extension of the presented results to the dynamical case in the spirit of  \cite{MFA}.}
}

\appendix
\subsection{Proof of Proposition \ref{prop_deltamax}}
For $\Sigma_D\in {\mathbf D}_n$, the Kullback-Leibler divergence  can be rewritten as
\begin{align*}&
2\mathcal{D}_{KL} (\Sigma_D \Vert \hat{\Sigma}) \nn\\
 =& \sum_{j=1}^{n}-\log d_j + \text{tr}(\text{diag}(d_1,...,d_n)\hat{\Sigma}^{-1}) + \log |\hat{\Sigma}|- n\\
= &  \left(\sum_{j=1}^{n}-\log d_j + d_j\gamma_j\right)+ \log |\hat{\Sigma}|- n.\\
\end{align*}
Thus the minimization problem in \eqref{min_diagonal} is equivalent to 
\begin{align*}
\min_{\Sigma_D\in {\mathbf D}_n} & 2\mathcal{D}_{KL} (\Sigma_D \Vert \hat{\Sigma})\nn\\
& =\left(\sum_{j=1}^{n} \min_{d_j} -\log d_j + d_j\gamma_j\right) + \log |\hat{\Sigma}| - n.
\end{align*}
Since $-\log d_j + d_j\gamma_j$ is convex with respect to $d_j$, by setting equal to zero the first derivative with respect to $d_j, j=1,...,n$, it easily follows that
\[\Sigma_D^{OPT}= \text{diag}(\gamma^{-1}_{1}, ..., \gamma^{-1}_{n}).\]
Then, $\delta_{max}$ can be determined as 
\begin{equation*}
\begin{aligned}
\delta_{max} & = 2\mathcal{D}_{KL}(\Sigma_D^{OPT} \Vert \hat{\Sigma})\\ 
& = \bigg{(} \sum_{j=1}^{n} -\log(\gamma^{-1}_{j}) +1  \bigg{)}  + \log |\hat{\Sigma}|- n \\
& = -\log|\text{diag}(\gamma^{-1}_{1}, ..., \gamma^{-1}_{n})|+ \log|\hat{\Sigma}|\\
& = \log | \dd(\hat{\Sigma}^{-1})\hat{\Sigma})|.
\end{aligned}
\end{equation*}

\subsection{Proof of Proposition \ref{prop_dual}}
By  substituting the obtained optimal conditions \eqref{der_var_wrt_r} and \eqref{der_var_wrt_sigma_opt} into \eqref{lagrangian}, we get the following expression where we have defined
$\Delta:=\lambda(\lambda\hat{\Sigma}^{-1}+ \ofd(\Theta)-\Gamma)^{-1}$:
\begin{align*}\label{dualR}
&J ( \lambda,\Gamma,\Theta) \\ &=\text{tr}(L)+ \lambda(-\log|\Delta|+\text{tr}(\hat{\Sigma}^{-1}\Delta)+\log|\hat{\Sigma}|-n-\delta)\\
&  \hspace{3mm} -\text{tr}((I+\Gamma-\ofd(\Theta))L)-\text{tr}(\Gamma\Delta)\\
& \hspace{3mm}+\text{tr}(\ofd(\Theta)\Delta)+\text{tr}(\Gamma L) -\text{tr}(\ofd(\Theta)L).
\end{align*}
That simplifies into:
\begin{equation*}
\label{dual-1}
\begin{aligned}
J(& \lambda, \Gamma,\Theta)\\ 
= \ & \lambda (-\log|(\hat{\Sigma}^{-1}+\lambda^{-1}( \ofd(\Theta)-\Gamma))^{-1}|+\log|\hat{\Sigma}|\\
& -n-\delta 
+ \text{tr}(\hat{\Sigma}^{-1}(\hat{\Sigma}^{-1}+\lambda^{-1}(\ofd(\Theta)-\Gamma))^{-1}))\\
& - \text{tr}(\Gamma(\hat{\Sigma}^{-1}+\lambda^{-1}( \ofd(\Theta)-\Gamma))^{-1})\\
&  + \text{tr}(\ofd(\Theta)(\hat{\Sigma}^{-1}+\lambda^{-1}(\ofd(\Theta)-\Gamma))^{-1})\\
= \ & \lambda(-\log|(\hat{\Sigma}^{-1}+\lambda^{-1}(\ofd(\Theta)-\Gamma))^{-1}|+ \log|\hat{\Sigma}|\\
& -n-\delta) 
+ \text{tr}(\lambda \hat{\Sigma}^{-1}(\hat{\Sigma}^{-1} +\lambda^{-1}( \ofd(\Theta)-\Gamma))^{-1})\\
& -\text{tr}(\Gamma(\hat{\Sigma}^{-1}+\lambda^{-1}(\ofd(\Theta)-\Gamma))^{-1})\\
& +\text{tr}(\ofd(\Theta)(\hat{\Sigma}^{-1}+\lambda^{-1}(\ofd(\Theta)-\Gamma))^{-1})\\
= \ & \lambda(-\log|(\hat{\Sigma}^{-1}+\lambda^{-1}(\ofd(\Theta)-\Gamma))^{-1}|+ \log|\hat{\Sigma}|\\
& -n -\delta) + \text{tr}((\lambda\hat{\Sigma}^{-1}+\ofd(\Theta)-\Gamma)\times \\
&(\hat{\Sigma}^{-1}+\lambda^{-1}(\ofd(\Theta)-\Gamma))^{-1})\\
= \ &  \lambda(-\log|(\hat{\Sigma}^{-1}+\lambda^{-1}(\ofd(\Theta)-\Gamma))^{-1}|+\log|\hat{\Sigma}|\\
& -n -\delta)
+\text{tr}(\lambda(\hat{\Sigma}^{-1}+\lambda^{-1}(\ofd(\Theta)-\Gamma))\times \\
& (\hat{\Sigma}^{-1}+\lambda^{-1}(\ofd(\Theta)-\Gamma))^{-1})\\
= \ &  \lambda(-\log|(\hat{\Sigma}^{-1}+\lambda^{-1}(\ofd(\Theta)-\Gamma))^{-1}|+\log|\hat{\Sigma}|\\
& -n -\delta)+n\lambda\\
= \ &  \lambda(\log|(\hat{\Sigma}^{-1}+\lambda^{-1}(\ofd(\Theta)-\Gamma))|+\log|\hat{\Sigma}|-\delta).
\end{aligned}
\end{equation*}

Since $J$ does not depend on $\Lambda$, we can eliminate it and, in view of (\ref{der_var_wrt_r}), condition $\Lambda\succeq 0$, is replaced by \alg{I+\Gamma-\ofd(\Theta)\succeq 0.}

\subsection{Proof of Proposition \ref{prop_hessian}}

Consider the function
\[\tilde{F}(\lambda,X)=-\lambda\log|\hat{\Sigma}+\lambda^{-1}X|.\]
Since $\tilde{F}(\lambda,X)$ differs from $F(\lambda,X)$ only by terms which are linear in $(\lambda,X)$ the second variations of the two functions are equivalent. Thus, in what follows we will focus on $\tilde{F}(\lambda,X)$.
The first variation of $\tilde{F}(\lambda,X)$ in direction $(\delta\lambda, \delta X)$ is 
\begin{align*}
\delta\tilde{F}& (\lambda,X;\delta \lambda, \delta X)=  -\log|\hat{\Sigma}+\lambda^{-1}X|\delta\lambda\\ &+\lambda^{-1}\text{tr}((\hat{\Sigma}+\lambda^{-1}X)^{-1}X)\delta\lambda - \text{tr}((\hat{\Sigma}+\lambda^{-1}X)^{-1}\delta X).
\end{align*}
The second variation of $\tilde{F}(\lambda,X)$ in direction $(\delta\lambda, \delta X)$ is 
\begin{align*}
\delta^{2}& \tilde{F}(\lambda, X;\delta \lambda, \delta X)= \\
& \lambda^{-1}\text{tr} \big{(}   (\hat{\Sigma}^{-1}+\lambda^{-1}X)^{-1} \delta X  (\hat{\Sigma}^{-1}+\lambda^{-1}X)^{-1}\delta X\big{)}\\
& -2\left[\lambda^{-2}\text{tr} \big{(}  (\hat{\Sigma}^{-1}+\lambda^{-1}X)^{-1}\delta X   (\hat{\Sigma}^{-1}+\lambda^{-1}X)^{-1}X \big{)}\delta\lambda\right]\\
& +\lambda^{-3}\text{tr} \big{(}  (\hat{\Sigma}^{-1}+\lambda^{-1}X)^{-1} X (\hat{\Sigma}^{-1}+\lambda^{-1}X)^{-1} X \big{)}\delta\lambda^{2}.
\end{align*}
Now, by using the Kronecker product and the vec operator and defining $x:=\text{vec}(X)$, $\delta x:=\text{vec}(\delta X)$, and $K:= (\hat{\Sigma}^{-1}+\lambda^{-1}X)^{-1}\otimes (\hat{\Sigma}^{-1}+\lambda^{-1}X)^{-1}$ the Hessian in Proposition \ref{prop_hessian} immediately follows.

\begin{IEEEbiography}
{Valentina Ciccone}  received the M.Sc. in Economics and Finance (cum laude) at the University of Padova, Italy, in 2016. She is currently a Ph.D student in Information Engineering at the University of Padova.
Her present research interests are in the areas of Systems Identification, Financial Econometrics and Stochastic Systems. 
\end{IEEEbiography}

\begin{IEEEbiography}
{Augusto Ferrante} was born in Piove di Sacco, Italy, on August 5, 1967. 
He received the ``Laurea" degree, {\em cum laude}, in Electrical Engineering in 1991
and the Ph.D. degree in Control Systems Engineering in 1995, both from the University of Padova.

He has been on the faculty of the Colleges of Engineering of the University of Udine,  and of the ``Politecnico di Milano". 
He is presently Professor in the ``Department of Information Engineering" of the University of Padova.

His research interests are in the areas of linear systems, spectral estimation, optimal control and optimal filtering, quantum control, and stochastic realization.
\end{IEEEbiography}

\begin{IEEEbiography}
{Mattia Zorzi} received the M.S. degree in Automation Engineering and the Ph.D. degree in Information Engineering from the University of Padova in 2009 and 2013, respectively. 
He held a postdoctoral position with the Department of Electrical Engineering and Computer Science, University of Liege (BE) in 2013-2014. He held a visiting position with the Department of Electrical and Computer Engineering, U.C. Davis (USA), and with the Department of Engineering, University of Cambridge (UK), in 2011 and 2013-2014, respectively. He is currently an Assistant Professor
with the Department of Information Engineering, University of Padova. Dr. Zorzi serves as a member of the Conference Editorial Board of the IEEE CSS. His current research interests include Machine Learning, Robust Estimation, Graphical Models and Identification Theory.\end{IEEEbiography}

\end{document}